\documentclass[11pt,a4]{article}
\usepackage{latexsym}
\usepackage{amssymb}
\usepackage{graphics}
\input epsf
\makeatletter
\def\section{\@startsection {section}{1}{\z@}{-3.5ex plus -1ex minus
 -.2ex}{2.3ex plus .2ex}{\large\sc}}
\def\subsection{\@startsection{subsection}{2}{\z@}{-3.25ex plus -1ex minus
 -.2ex}{1.5ex plus .2ex}{\normalsize\sc}}
\makeatother

\makeatletter \@addtoreset{equation}{section}

\makeatother

\newcommand{\nc}{\newcommand}
\newcommand{\rnc}{\renewcommand}

\nc{\bea}{\begin{eqnarray}} \nc{\eea}{\end{eqnarray}}
\nc{\be}{\bea} \nc{\ee}{\eea}

\rnc{\a}{\alpha} \nc{\ab}{\bar{\a}} \nc{\ap}{\a^{+}}
\nc{\abm}{\ab^{-}} \rnc{\b}{\beta} \nc{\bb}{\bar{\b}}
\nc{\bbp}{\bb_{\zb}^{+}} \nc{\bm}{\b_{z}^{-}}
\nc{\oa}{\overline{\a}} \nc{\ob}{\overline{\b}} \rnc{\gg}{\gamma}
\rnc{\d}{\delta} \nc{\f}{\phi} \nc{\fb}{\bar{\phi}}
\nc{\vf}{\varphi} \nc{\p}{\psi}

\rnc{\c}{\chi} \nc{\la}{\lambda} \nc{\m}{\mu} \nc{\n}{\nu}
\rnc{\o}{\omega} \nc{\Om}{\Omega} \rnc{\t}{\theta}
\nc{\eps}{\epsilon} \rnc{\S}{\Sigma} \nc{\F}{\Phi}

\nc{\trac}[2]{{\textstyle\frac{#1}{#2}}}

\nc{\ex}[1]{\mbox{e}^{\,\textstyle#1}}

\nc{\som}[9]{(\begin{array}{ccc}#1&#2&#3\\#4&#5&#6\\#7&#8&#9%
\end{ar\nc{\mat}[4]{(\begin{array}{cc}#1&#2\\#3&#4\end{array})}
ray})}

\nc{\tr}{\mathop{\mbox{tr}}\nolimits}
\nc{\ad}{\mathop{\mbox{ad}}\nolimits}
\nc{\Tr}{\mathop{\mbox{Tr}}\nolimits}
\nc{\Det}{\mathop{\mbox{Det}}\nolimits}
\nc{\rk}{\mathop{\mbox{rk}}\nolimits} \nc{\ra}{\rightarrow}
\nc{\Ra}{\Rightarrow} \nc{\LRa}{\Leftrightarrow} \nc{\ot}{\otimes}
\rnc{\ss}{\subset} \nc{\nul}{\noindent\underline}
\nc{\non}{\nonumber\\} \nc{\ZZ}{\mathbb{Z}} \nc{\RR}{\mathbb{R}}
\nc{\CC}{\mathbb{C}}

\nc{\subs}[1]{{\vspace*{0.5cm}}%
{\noindent\underline{#1}}%
{\vspace*{0.3cm}}}

\nc{\zb}{\bar{z}}
\rnc{\lg}{\mathfrak{g}}
\nc{\lt}{\mathfrak{t}}
\nc{\lk}{\mathfrak{k}}
\nc{\lh}{\mathfrak{h}} \nc{\pik}{\Pi_{\lk}} \nc{\pip}{\Pi_{+}}
\nc{\pim}{\Pi_{-}} \nc{\pih}{\Pi_{\lh}} \nc{\jz}{J_{z}}
\nc{\jzh}{\jz^{\lh}} \nc{\jzp}{\jz^{+}} \nc{\jzm}{\jz^{-}}
\nc{\del}{\partial} \nc{\dz}{\del_{z}} \nc{\dzb}{\del_{\bar{z}}}
\nc{\az}{A_{z}} \nc{\azb}{A_{\bar{z}}} \nc{\g}{g^{-1}}
\nc{\dw}{\Delta_{W}} \nc{\Ad}{{\mbox{Ad}}}
\nc{\ks}{Ka\-za\-ma-\-Su\-zu\-ki} \nc{\KS}{\ks} \nc{\ksm}{\ks\
model} \rnc{\AA}{{\Bbb A}} \nc{\BB}{{\Bbb B}} \nc{\PP}{{\Bbb P}}
\nc{\cpm}{\CC\PP(m)} \nc{\cpn}{\CC\PP(n)} \nc{\cp}[1]{\CC\PP(#1)}
\nc{\gmn}{G(m,m+n)} \nc{\gmnk}{\gmn_{k}} \nc{\cO}{{\cal O}}
\nc{\bcO}{\bar{\cO}} \nc{\bO}{\bar{O}} \nc{\oQ}{\overline{Q}}

\newtheorem{theorem}{Theorem}

\newtheorem{lemma}{Lemma}
\newtheorem{corollary}{Corollary}

\newtheorem{example}{Example}
\newtheorem{remark}{Remark}

\begin{document}
\thispagestyle{empty}
\title{Gaps in Nonsymmetric Numerical Semigroups}
\author{Leonid G. Fel ${}^{\dag}$ and Francesca Aicardi ${}^{\ddag}$
\\ \\ \\
${}^{\dag}$ Department of Civil Engineering, Technion, Haifa
32000, Israel\\ {\sl e-mail: lfel@tx.technion.ac.il}\\ and \\
${}^{\ddag}$ Sistiana 56 PR, Trieste, Italy\\{\sl e-mail:
aicardi@sissa.it}}
\date{\today}

\maketitle
\begin{abstract}
There exist two different sorts of  gaps in the nonsymmetric
numerical additive semigroups  finitely generated by a minimal set
of positive integers $\{d_1,\ldots,d_m\}$. The $h$--gaps are specific
only for the nonsymmetric semigroups while the $g$--gaps are possessed by
both, symmetric and nonsymmetric semigroups. We derive the generating
functions for the corresponding sets of gaps, $\Delta_{{\cal H}}\left(
{\bf d}^m\right)$ and $\Delta_{{\cal G}}\left({\bf d}^m\right)$, and prove
several statements on the minimal and maximal values of the $h$--gaps.
Detailed description of both sorts of gaps is given for three special
kinds of nonsymmetric semigroups: semigroups with maximal embedding
dimension, semigroups of maximal and almost maximal length, and
pseudo--symmetric semigroups.
\end{abstract}
\section{Introduction}\label{p1}
In this paper we study the additive numerical semigroup in ${\mathbb N}
\cup\{0\}$ generated by a finite set of positive integers $\left\{d_1,
\ldots,d_m\right\}$ and the arrangement of the set of the integers which
are unrepresentable by $\left\{d_1,\ldots,d_m\right\}$. Such integers
are also known as the {\em  gaps} of the numerical semigroup.

First, recall the main known facts \cite{baru97}. Semigroup ${\sf S}(d_1,
\ldots,d_m)$,
\begin{equation}
{\sf S}(d_1,\ldots,d_m)=\left\{s\in {\mathbb N}\cup\{0\}\quad |\quad
s=\sum_{i=1}^m x_i d_i,\quad x_i\in {\mathbb N}\cup\{0\}\right\}\;,\label{tr1}
\end{equation}
is said to be generated by the minimal set of $m$ natural numbers $d_1<\ldots
<d_m$, if there are no nonnegative integers $b_{i,j}$ for which the following
linear dependence holds :
\begin{equation}
d_i=\sum_{j\neq i}^mb_{i,j}d_j\;,\;\;\;b_{i,j}\in\{0,1,\ldots\}\;\;\;
\mbox{for any}\;\;i\leq m\;,\label{deff1}
\end{equation}
and $\gcd(d_1,\ldots,d_m)=1$. It is classically known that $d_1\geq m$
\cite{abhy67}. For short we denote the tuple $(d_1,\ldots,d_m)$ by ${\bf d}
^m$ where $m$ is called {\em the embedding dimension} of the semigroup.
Henceforth ${\bf d}^m$ will be a minimal generating set.

The least positive integer ($d_1$) belonging to ${\sf S}\left({\bf d}^m\right)$
is called {\em the multiplicity}. The {\em  conductor} $c\left({\bf d}^m
\right)$ of ${\sf S}\left({\bf d}^m\right)$ is defined by $c\left({\bf
d}^m\right):=\min\left\{s\in {\sf S}\left({\bf d}^m\right)\;|\;s+{\mathbb N}
\cup\{0\}\subset {\sf S}\left({\bf d}^m\right)\right\}$. Denote by $\Delta
\left({\bf d}^m\right)$ the complement of ${\sf S}\left({\bf d}^m\right)$
in ${\mathbb N}$, i.e. $\Delta\left({\bf d}^m\right)={\mathbb N}\setminus
{\sf S}\left({\bf d}^m\right)$. It is the set of gaps. The cardinality
($\#$) of $\Delta\left({\bf d}^m\right)$ is called the {\em genus}
of ${\sf S}\left({\bf d}^m\right)$,
\begin{eqnarray}
G\left({\bf d}^m\right):=\#\Delta\left({\bf d}^m\right)\;.\label{tr2}
\end{eqnarray}
Introduce the generating functions $\Phi\left({\bf d}^m;z\right)$ for the
set $\Delta\left({\bf d}^m\right)$ of unrepresentable integers and $H\left({
\bf d}^m;z\right)$ for the set of representable integers
\begin{eqnarray}
\Phi\left({\bf d}^m;z\right)=\sum_{s\;\in\;\Delta\left({\bf d}^m\right)}z^s\;,
\;\;\;\;\;\;\;\;H\left({\bf d}^m;z\right)=\sum_{s\;\in\;{\sf S}\left({\bf
d}^m\right)}z^s\;.\label{tr3}
\end{eqnarray}
The latter is referred to as {\em the Hilbert series}. The function $\Phi
\left({\bf d}^m;z\right)$ determines the Frobenius number and the genus
as follows :
\begin{eqnarray}
F\left({\bf d}^m\right)=c\left({\bf d}^m\right)-1=\deg\Phi\left({\bf d}^m;z
\right)\;,\;\;\;G\left({\bf d}^m\right)=\Phi\left({\bf d}^m;1\right)\;.
\label{tr4}
\end{eqnarray}
In terms of the rings theory the entities $F\left({\bf d}^m\right)$ and $G
\left({\bf d}^m\right)$ are also referred to as {\em the Castelnuovo--Mumford
regularity} \cite{stur05} and {\em the length} \cite{baru97} of semigroups,
respectively.

A semigroup ${\sf S}\left({\bf d}^m\right)$ is called {\em symmetric} iff
for any integer $s$ the following condition holds :
\begin{eqnarray}
s\in {\sf S}\left({\bf d}^m\right)\;\;\;\Longleftrightarrow\;\;\;F\left({\bf
d}^m\right)-s\not\in {\sf S}\left({\bf d}^m\right)\;.\label{tr5}
\end{eqnarray}
Otherwise ${\sf S}\left({\bf d}^m\right)$ is called {\em nonsymmetric}. $F
\left({\bf d}^m\right)$ and $G\left({\bf d}^m\right)$ are related as follows :
\begin{eqnarray}
2G\left({\bf d}^m\right)=c\left({\bf d}^m\right)\;\;\mbox{if}\;\;{\sf S}
\left({\bf d}^m\right)\;\;\mbox{is symmetric, and}\;\;2G\left({\bf d}^m\right)>
c\left({\bf d}^m\right)\;\;\mbox{ otherwise}.\label{tr6}
\end{eqnarray}
Denote by $t\left({\bf d}^m\right)$ {\em the type} of the numerical semigroup
${\sf S}\left({\bf d}^m\right)$ which coincides with the cardinality of set
${\sf S}^{\prime}\left({\bf d}^m\right)$ \cite{frob87}, where
\begin{eqnarray}
{\sf S}^{\prime}\left({\bf d}^m\right)=\left\{F_i\in {\mathbb Z}\;|\;F_i\not\in
{\sf S}\left({\bf d}^m\right),\;F_i+s\in {\sf S}\left({\bf d}^m\right),\;
s\in {\sf S}\left({\bf d}^m\right)\setminus \{0\},\;i=1,\ldots t\left({\bf d}
^m\right)\right\}.\;\;\label{tr6a}
\end{eqnarray}
Set ${\sf S}^{\prime}\left({\bf d}^m\right)$ is not empty since $F\left(
{\bf d}^m\right)\in {\sf S}^{\prime}\left({\bf d}^m\right)$ for any
minimal generating set $(d_1,\ldots ,d_m)$. Due to \cite{frob87},
Proposition 2, semigroup ${\sf S}\left({\bf d}^m\right)$ is symmetric iff
$t\left({\bf d}^m\right)=1$.
\section{The Two Sorts of Gaps in Numerical Semigroups}\label{p2}
Due to (\ref{tr5}) set $\Delta\left({\bf d}^m\right)$ can be decomposed
into two sorts,
\begin{eqnarray}
\Delta\left({\bf d}^m\right)=\Delta_{{\cal G}}\left({\bf d}^m\right)\cup
\Delta_{{\cal H}}\left({\bf d}^m\right)\;,\;\;\;\Delta_{{\cal G}}\left(
{\bf d}^m\right)\cap\Delta_{{\cal H}}\left({\bf d}^m\right)=\emptyset\;,
\label{tr7}
\end{eqnarray}
where the sets of the $g$--gaps and the $h$--gaps are defined, respectively, by
\begin{eqnarray}
&&\Delta_{{\cal G}}\left({\bf d}^m\right)=\left\{g\not\in {\sf S}\left({\bf
d}^m\right)\;|\;F\left({\bf d}^m\right)-g\in {\sf S}\left({\bf d}^m\right)
\right\},\;\;\;\#\Delta_{{\cal G}}\left({\bf d}^m\right)=c\left({\bf
d}^m\right)-G\left({\bf d}^m\right),\;\;\;\;\;\;\;\;\;\;\;\;\label{tr8}\\
&&\Delta_{{\cal H}}\left({\bf d}^m\right)=\left\{h\not\in {\sf S}\left({\bf
d}^m\right)\;|\;F\left({\bf d}^m\right)-h\not\in {\sf S}\left({\bf d}^m\right)
\right\},\;\;\;\#\Delta_{{\cal H}}\left({\bf d}^m\right)=2G\left({\bf d}^m
\right)-c\left({\bf d}^m\right).\;\;\;\;\;\;\;\;\;\;\;\;\label{tr9}
\end{eqnarray}
\centerline{\epsfbox{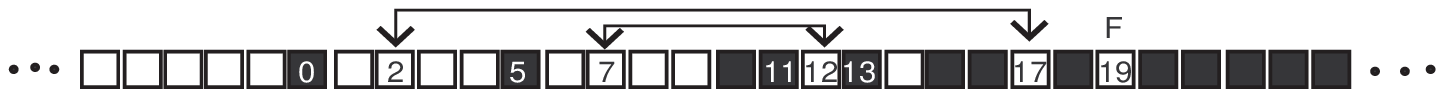}} {\small Figure.1:  Initial part
(black) of the semigroup generated by (5,11,13). The arrows show
the pairs of $h$--gaps. More details about this semigroup will be
given later, in Example \ref{ex1}.}

Following \cite{stan96} recall the relation  between $\Phi\left({\bf d}^m;z
\right)$ and the Hilbert series $H\left({\bf d}^m;z\right)$ of the graded
subring ${\sf k}\left[z^{d_1},\ldots,z^{d_m}\right]$ associated with
semigroup ${\sf S}\left({\bf d}^m\right)$ :
\begin{eqnarray}
\Phi\left({\bf d}^m;z\right)=\frac{1}{1-z}-H\left({\bf d}^m;z\right)\;,\;\;
\;\;\;H\left({\bf d}^m;z\right)=\frac{Q\left({\bf d}^m;z\right)}{\prod_{j=1}^m
\left(1-z^{d_j}\right)}\;,\label{tr16}
\end{eqnarray}
where $H\left({\bf d}^m;z\right)$ has a pole $z=1$ of order 1. The numerator
$Q\left({\bf d}^m;z\right)$ is a polynomial in $z$ of the form
\begin{eqnarray}
Q\left({\bf d}^m;z\right)=1-\sum_{j=1}^{\beta_1\left({\bf d}^m\right)}z^{C_{j,
1}}+\sum_{j=1}^{\beta_2\left({\bf d}^m\right)}z^{C_{j,2}}-\ldots\pm\sum_{j=1}^
{\beta_{m-1}\left({\bf d}^m\right)}z^{C_{j,m-1}}\;.\label{bet1}
\end{eqnarray}
The summands of the powers $z^{C_{j,i}}$ in the last formula stand for the
syzygies of different kinds and $C_{j,i}$, $C_{j+1,i}>C_{j,i}$, are the
degrees of homogeneous basic invariants for the syzygies of the $i$th kind.
The numbers of the terms $z^{C_{j,i}}$ in the summands are determined by the
Betti numbers $\beta_i\left({\bf d}^m\right)$ which satisfy the relation
\begin{eqnarray}
1-\beta_1\left({\bf d}^m\right)+\beta_2\left({\bf d}^m\right)-\ldots
\pm\beta_{m-1}\left({\bf d}^m\right)=0\;.\label{bet2}
\end{eqnarray}
It is also known \cite{heku71}, \cite{brown92} that
\begin{eqnarray}
\beta_{m-1}\left({\bf d}^m\right)=t\left({\bf d}^m\right)\;.\label{bet2n}
\end{eqnarray}
\begin{theorem}\label{the1}{\rm $\;$}\\
The generating functions for the sets $\Delta_{{\cal H}}\left({\bf d}^m
\right)$ and $\Delta_{{\cal G}}\left({\bf d}^m\right)$ are given by
\begin{eqnarray}
\sum_{h\;\in\;\Delta_{{\cal H}}\left({\bf d}^m\right)}z^h&=&-H\left({\bf d}^m;
z\right)-H\left({\bf d}^m;\frac{1}{z}\right)\cdot z^{F\left({\bf d}^m\right)}\;,
\label{tr17a}\\
\sum_{g\;\in\;\Delta_{{\cal G}}\left({\bf d}^m\right)}z^g&=&\frac{1}{1-z}+
H\left({\bf d}^m;\frac{1}{z}\right)\cdot z^{F\left({\bf d}^m\right)}\;.
\label{tr17b}
\end{eqnarray}
\end{theorem}
{\sf Proof} $\;\;\;$We have by (\ref{tr7})
\begin{eqnarray}
\Phi\left({\bf d}^m;z\right)=\sum_{g\;\in\;\Delta_{{\cal G}}\left({\bf
d}^m\right)}z^g+\sum_{h\;\in\;\Delta_{{\cal H}}\left({\bf d}^m\right)}z^h\;.
\label{tr10}
\end{eqnarray}
Consider the transformation
\begin{eqnarray}
\Phi\left({\bf d}^m;\frac{1}{z}\right)\cdot z^{F({\bf d}^m)}=\sum_{g\;\in\;
\Delta_{{\cal G}}({\bf d}^m)}z^{F({\bf d}^m)-g}+\sum_{h\;\in\;\Delta_{{\cal
H}}({\bf d}^m)}z^{F({\bf d}^m)-h}\;.\label{tr11}
\end{eqnarray}
However, according to (\ref{tr9}),
\begin{eqnarray}
F\left({\bf d}^m\right)-g\in {\sf S}\left({\bf d}^m\right)\;,\;\;\;
F\left({\bf d}^m\right)-h\in\Delta_{{\cal H}}\left({\bf d}^m\right)\;.
\label{tr12}
\end{eqnarray}
Thence,
\begin{eqnarray}
\Phi\left({\bf d}^m;\frac{1}{z}\right)\cdot z^{F({\bf d}^m)}=\sum_{
s\;\in\;[0;F({\bf d}^m)]}z^s+\sum_{h\;\in\;\Delta_{{\cal H}}
({\bf d}^m)}z^h\;.\label{tr13}
\end{eqnarray}
\begin{remark}\label{re0}
{\rm The elements in ${\sf S}\left({\bf d}^m\right)$ which are less than
$F\left({\bf d}^m\right)$, i.e. $s\in {\sf S}\left({\bf d}^m\right)\cap
[0;F\left({\bf d}^m\right)]$, are also known as {\em the nongaps} of the
semigroup. Their cardinality coincides with $\#\Delta_{{\cal G}}\left({\bf
d}^m\right)$.}
\end{remark}
Making summation of (\ref{tr10}) and (\ref{tr13}), we get
\begin{eqnarray}
\Phi\left({\bf d}^m;z\right)+\Phi\left({\bf d}^m;\frac{1}{z}\right)\cdot z^{F
\left({\bf d}^m\right)}=\sum_{s\;\in\;[0;F({\bf d}^m)]}z^s+\sum_{g
\;\in\;\Delta_{{\cal G}}({\bf d}^m)}z^{g}+2\sum_{h\;\in\;
\Delta_{{\cal H}}({\bf d}^m)}z^{h}\;.\label{tr14}
\end{eqnarray}
So, we come to the generating functions of the set $\Delta_{{\cal H}}\left(
{\bf d}^m\right)$,
\begin{eqnarray}
\sum_{h\;\in\;\Delta_{{\cal H}}({\bf d}^m)}z^h=\Phi\left({\bf d}^m;z\right)+
\Phi\left({\bf d}^m;\frac{1}{z}\right)\cdot z^{F({\bf d}^m)}-
\sum_{k=0}^{F({\bf d}^m)}z^k\;.\label{tr15}
\end{eqnarray}
Substituting the 1st equality of (\ref{tr16}) into (\ref{tr15}) we obtain
\begin{eqnarray}
\sum_{h\;\in\;\Delta_{{\cal H}}({\bf d}^m)}z^h&=&\frac{1}{1-z}-
\frac{z^{F({\bf d}^m)+1}}{1-z}-\sum_{k=0}^{F({\bf d}^m)}z^k-
H\left({\bf d}^m;z\right)-H\left({\bf d}^m;\frac{1}{z}\right)\cdot z^{F
({\bf d}^m)}\nonumber\\
&=&-H\left({\bf d}^m;z\right)-H\left({\bf d}^m;\frac{1}{z}\right)\cdot
z^{F({\bf d}^m)}\;.\label{tr17}
\end{eqnarray}
Finally, combining (\ref{tr16}), (\ref{tr10}) and (\ref{tr17a}) we get formula
(\ref{tr17b}).$\;\;\;\;\;\;\Box$

Theorem \ref{the1} immediately results in the following duality relation.
\begin{corollary}\label{cor01}{\rm $\;$}\\
The semigroup ${\sf S}\left({\bf d}^m\right)$ is symmetric iff
\begin{eqnarray}
H\left({\bf d}^m;z\right)+H\left({\bf d}^m;\frac{1}{z}\right)\cdot z^{F({\bf 
d}^m)}=0\;.\label{tr17f}
\end{eqnarray}
\end{corollary}
\vspace{.3cm}
Equations (\ref{tr17a}) and (\ref{tr17b}) allow us to formulate two new
theorems. Let a semigroup ${\sf S}\left({\bf d}^m\right)$ be given, and the
Hilbert series $H\left({\bf d}^m;z\right)$ of its graded subring ${\sf k}
\left[z^{d_1},\ldots, z^{d_m}\right]$ be given by (\ref{tr16}) and
(\ref{bet1}). Define polynomial $V\left({\bf d}^m;z\right)$ by
\begin{eqnarray}
V\left({\bf d}^m;z\right)=(-1)^{m-1}z^{\deg Q({\bf d}^m;z)}Q\left({\bf d}
^m;\frac{1}{z}\right)-Q\left({\bf d}^m;z\right)\;.\label{tr18}
\end{eqnarray}
\begin{theorem}\label{the2}{\rm $\;$}\\
The semigroup ${\sf S}\left({\bf d}^m\right)$ is symmetric iff the polynomial
$V\left({\bf d}^m;z\right)$ is vanishing and ${\sf S}\left({\bf d}^m\right)$
is nonsymmetric iff $V\left({\bf d}^m;z\right)$ is divided by $\prod_{j=1}^
m\left(1-z^{d_j}\right)$.
\end{theorem}
{\sf Proof} $\;\;\;$ Combining (\ref{tr18}) with (\ref{tr17a}) and
(\ref{tr16}), we obtain:
\begin{eqnarray}
V\left({\bf d}^m;z\right)=\prod_{j=1}^m\left(1-z^{d_j}\right)\sum_{h\;\in\;
\Delta_{{\cal H}}\left({\bf d}^m\right)}z^h\;.\label{tr19a}
\end{eqnarray}
The statement of the Theorem follows immediately from this identity.
$\;\;\;\;\;\;\Box$

Denote by $\min\Delta_{{\cal H}}\left({\bf d}^m\right)$ and $\max\Delta_{{\cal
H}}\left({\bf d}^m\right)$ the least and the largest $h$--gaps of semigroup 
${\sf S}\left({\bf d}^m\right)$, respectively.
\begin{theorem}\label{the3}{\rm $\;$}\\
Let a nonsymmetric semigroup ${\sf S}\left({\bf d}^m\right)$ be given, and
the Hilbert series $H\left({\bf d}^m;z\right)$ of its graded subring be given
by (\ref{tr16}), (\ref{bet1}). Then $\min\Delta_{{\cal H}}\left({\bf d}^m
\right)$ is given by the minimum degree among the terms of polynomial $V\left(
{\bf d}^m;z\right)$, and $\max\Delta_{{\cal H}}\left({\bf d}^m\right)$ is
given by
\begin{equation}
\max\Delta_{{\cal H}}\left({\bf d}^m\right)=\deg V\left({\bf d}^m;z\right)
-\Sigma_m\;,\label{tr22}
\end{equation}
where $\Sigma_m=\sum_{j=1}^md_j$.
\end{theorem}
{\sf Proof} $\;\;\;$
The term with minimum degree of the sum $\sum_{h\;\in\;\Delta_{{\cal H}}
\left({\bf d}^m\right)}z^h$ has to coincide with the corresponding term in
the left hand side of (\ref{tr19a}). Thus, the minimal $h$--gap comes from the
term with minimal degree in $V\left({\bf d}^m;z\right)$. On the other hand,
the maximal degree among the terms of $V\left({\bf d}^m;z\right)$ satisfies:
\[\deg V\left({\bf d}^m;z\right)=\max\Delta_{{\cal H}}\left({\bf d}^m\right)
+\Sigma_m\;,\]
that concludes the proof.$\;\;\;\;\;\;\Box$
\begin{remark}\label{rem00}
{\rm The term with degree  $\min\Delta_{{\cal H}}\left({\bf d}^m\right)$ has
coefficient equal to $+1$, whereas the term with degree $\max\Delta_{{\cal H}}
\left({\bf d}^m\right)+\Sigma_m$ has coefficient equal to $(-1)^{m+1}$.}
\end{remark}
\begin{corollary}
\begin{equation}
\min\Delta_{{\cal H}}\left({\bf d}^m\right)\le\#\Delta_{{\cal G}}\left({\bf
d}^m\right)\;,\label{fr12}
\end{equation}
\end{corollary}
{\sf Proof} $\;\;\;$Since all the $h$--gaps are inside $\Delta_{{\cal H}}$, 
its cardinality satisfies
\begin{equation}
\#\Delta_{{\cal H}}\left({\bf d}^m\right)\le\max\Delta_{{\cal H}}\left({\bf
d}^m\right)-\min\Delta_{{\cal H}}\left({\bf d}^m\right)+1\;.\label{fr10}
\end{equation}
By the definition (\ref{tr9}) we obtain
\begin{equation}
\max\Delta_{{\cal H}}\left({\bf d}^m\right)=F\left({\bf d}^m\right)-
\min\Delta_{{\cal H}}\left({\bf d}^m\right)\;,\label{minmax}
\end{equation}
and, by the second equality in (\ref{tr9}), inequality (\ref{fr10}) reads:
\begin{equation}
2G\left({\bf d}^m\right)-c\left({\bf d}^m\right)\le F\left({\bf d}^m\right)-
\min\Delta_{{\cal H}}\left({\bf d}^m\right)-\min\Delta_{{\cal H}}\left({\bf
d}^m\right)+1=c\left({\bf d}^m\right)-2\min\Delta_{{\cal H}}\left({\bf
d}^m\right).\label{fr11}
\end{equation}
Taking into account (\ref{tr8}) we get finally (\ref{fr12}).
$\;\;\;\;\;\;\Box$

We finish this Section by discussing a question which is related to Theorem 
\ref{the1} : does it exist a symmetric semigroup ${\sf S}\left({\bf p}^n
\right)$ with minimal generating set ${\bf p}^n$ such that its total set of 
gaps $\Delta\left({\bf p}^n\right)$ coincides with the set of $g$--gaps 
$\Delta_{{\cal G}}\left({\bf d}^m\right)$ of a given nonsymmetric semigroup 
${\sf S}\left({\bf d}^m\right)$ ?
\begin{theorem}\label{the3a}{\rm $\;$}\\
Let a nonsymmetric semigroup ${\sf S}\left({\bf d}^m\right)$ be given with 
Hilbert series $H\left({\bf d}^m;z\right)$. It does not exist any symmetric 
semigroup ${\sf S}\left({\bf p}^n\right)$ with minimal generating set ${\bf 
p}^n$  such that $\Delta\left({\bf p}^n\right)=\Delta_{{\cal G}}\left({\bf 
d}^m\right)$.
\end{theorem}
{\sf Proof} $\;\;\;$Let, by way of contradiction, such symmetric semigroup
${\sf S}\left({\bf p}^n\right)$ does exist, i.e. $\Delta\left({\bf p}^n\right)=
\Delta_{{\cal G}}\left({\bf d}^m\right)$. In accordance with (\ref{tr3}),
(\ref{tr16}) and Theorem \ref{the1} we have
\begin{eqnarray}
\frac{1}{1-z}+H\left({\bf d}^m;\frac{1}{z}\right)\cdot z^{F({\bf d}^m)}=
\frac{1}{1-z}-H\left({\bf p}^n;z\right)\;.\label{add1}
\end{eqnarray}
Being symmetric, semigroup ${\sf S}\left({\bf p}^n\right)$ satisfies Corollary 
\ref{cor01}, and therefore (\ref{add1}) can be written as
\begin{eqnarray}
H\left({\bf d}^m;\frac{1}{z}\right)\cdot z^{F({\bf d}^m)}=H\left({\bf p}^n;
\frac{1}{z}\right)\cdot z^{F({\bf p}^n)}\;,\label{add2}
\end{eqnarray}
which requires necessarily $F\left({\bf d}^m\right)=F\left({\bf p}^n\right)$
since both Hilbert series start with unity. Inverting $z\to z^{-1}$ in the 
remained equality (\ref{add2}), $H\left({\bf d}^m;z^{-1}\right)=H\left({\bf 
p}^n;z^{-1}\right)$, and substituting the result into (\ref{add1}), we get
\begin{eqnarray}
H\left({\bf d}^m;z\right)+H\left({\bf d}^m;\frac{1}{z}\right)\cdot
z^{F({\bf d}^m)}=0\;.\label{add3}
\end{eqnarray}
According to Corollary \ref{cor01} we can conclude that semigroup ${\sf S}
\left({\bf d}^m\right)$ is symmetric. However, this contradicts our assumption 
about ${\sf S}\left({\bf d}^m\right)$ and proves the theorem.
$\;\;\;\;\;\;\Box$
\section{Low -- Dimensional Cases}\label{p3}
In this Section we apply Theorem \ref{the1} and Theorem \ref{the3} to 
numerical semigroups generated by two and three positive integers, respectively.
\subsection{Semigroup ${\sf S}\left({\bf d}^2\right)$}\label{s1p3}
Semigroup ${\sf S}\left({\bf d}^2\right)$ is always symmetric \cite{aper46},
and its Hilbert series reads
\begin{eqnarray}
H\left({\bf d}^2;z\right)=\frac{1-z^{d_1d_2}}{\left(1-z^{d_1}\right)\left(
1-z^{d_2}\right)}\;.\label{tr23}
\end{eqnarray}
Simple calculation of (\ref{tr17a}) in accordance with Theorem \ref{the2}
yields,
\begin{eqnarray}
\sum_{h\;\in\;\Delta_{{\cal H}}\left({\bf d}^2\right)}z^h=0\;.\label{tr118}
\end{eqnarray}
\subsection{Nonsymmetric semigroup ${\sf S}\left({\bf d}^3\right)$}
\label{s2p3a}
Nonsymmetric semigroups ${\sf S}\left({\bf d}^3\right)$ were studied by
algebraic means in \cite{herz70}, \cite{frob94}, \cite{denh03} and \cite{fel04}.
Recall its main results following \cite{fel04}.

Let ${\sf S}(d_1,d_2,d_3)\subset {\mathbb Z}_+\cup\{0\}$ be the
additive nonsymmetric numerical semigroup finitely generated by the minimal
set of positive integers $d_1<d_2<d_3$. Following Johnson \cite{john60}
define {\em the minimal relation} for the triple ${\bf d}^3=(d_1,d_2,d_3)$ :
\begin{eqnarray}
&&a_{11}d_1=a_{12}d_2+a_{13}d_3\;,\;\;\;a_{22}d_2=a_{21}d_1+a_{23}d_3\;,\;\;
\;a_{33}d_3=a_{31}d_1+a_{32}d_2\;,\label{joh1}
\end{eqnarray}
where
\begin{eqnarray}
&&a_{jj}=\min\left\{v_{jj}\;|\;v_{jj}\geq 2,\;v_{jj}d_j=v_{jk}d_k+
v_{jl}d_l,\;v_{jk},v_{jl}\in {\mathbb Z}_+\right\}\;,\label{joh1a}\\
&&\gcd(a_{jj},a_{jk},a_{jl})=1\;,\;\;\;\mbox{and}\;\;\;(j,k,l)=(1,2,3),\;
(2,3,1),\;(3,1,2)\;.\nonumber
\end{eqnarray}
The uniquely defined values of $v_{ij},i\neq j$ which give $a_{ii}$ will be
denoted by $a_{ij},\ i\neq j$. The degeneracy of the matrix $(a_{ij})$
together with (\ref{joh1a}) results in strong equalities relating the matrix
elements $a_{ij}$ and the generators $d_k$: for any permutation of indices
$(i,j,k)$, $i,j,k=1,2,3$  the following identities hold \cite{fel04},
\cite{john60} :
\begin{eqnarray}\label{eq2}
a_{ii}=a_{ji}+a_{ki}\;,\;\;\;\;\;\;\;a_{ii}a_{jj}=d_k +a_{ij}a_{ji}\;.
\label{fr1}
\end{eqnarray}
The Hilbert series $H\left({\bf d}^3;z\right)$, the type $t\left({\bf d}^3
\right)$, the Frobenius number $F\left({\bf d}^3;z\right)$ and the genus
$G\left({\bf d}^3\right)$ read \cite{fel04} :
\begin{eqnarray}
H\left({\bf d}^3;z\right)=\frac{Q\left({\bf d}^3;z\right)}{\left(1-z^{d_1}
\right)\left(1-z^{d_2}\right)\left(1-z^{d_3}\right)}\;,\label{tr24}
\end{eqnarray}
where
\begin{eqnarray}
Q\left({\bf d}^3;z\right)&=&1-\sum_{i=1}^3z^{a_{ii}d_i}+z^{1/2\left[\langle{
\bf a},{\bf d}\rangle-J\left({\bf d}^3\right)\right]}+z^{1/2\left[\langle{
\bf a},{\bf d}\rangle+J\left({\bf d}^3\right)\right]}\;,\;\;\;\;t\left({\bf
d}^3\right)=2\;,\label{tr24a}\\
F\left({\bf d}^3\right)&=&\frac{1}{2}\left[\langle{\bf a},{\bf d}\rangle+
J\left({\bf d}^3\right)\right]-\Sigma_3\;,\;\;\;\;
G\left({\bf d}^3\right)=\frac{1}{2}\left(1+\langle{\bf a},{\bf d}\rangle
-\prod_{i=1}^3a_{ii}-\Sigma_3\right),\;\;\;\;\;\;\;\label{tr24b}\\
J^2\left({\bf d}^3\right)&=&\langle{\bf a},{\bf d}\rangle^2-4
\sum_{i>j}^3a_{ii}a_{jj}d_id_j+4d_1d_2d_3\;,\;\;\;\;
\langle{\bf a},{\bf d}\rangle=\sum_{i=1}^3a_{ii}d_i\;.\label{tr24c}
\end{eqnarray}
We prove the following Theorem.
\begin{theorem}\label{the4}{\rm $\;$}\\
The minimal and the maximal values in $\Delta_{{\cal H}}\left({\bf d}^3\right)$
are given by
\begin{eqnarray}
\min\Delta_{{\cal H}}\left({\bf d}^3\right)=J\left({\bf d}^3\right)\;,\;\;\;
\;\;\max\Delta_{{\cal H}}\left({\bf d}^3\right)=F\left({\bf d}^3\right)-
J\left({\bf d}^3\right)\;.\label{tr25}
\end{eqnarray}
\end{theorem}
{\sf Proof} $\;\;\;$Write the numerator $Q\left({\bf d}^3;z\right)$ in the
following form :
\begin{eqnarray}
Q\left({\bf d}^3;z\right)=1-\sum_{i+1}^3 z^{a_{ii}d_i}+z^{F({\bf d}^3)-
J\left({\bf d}^3\right)+\Sigma_3}+z^{F\left({\bf d}^3\right)+\Sigma_3}\;,
\label{ar2}
\end{eqnarray}
that in accordance with (\ref{tr18}) gives polynomial $V\left({\bf d}^3;z
\right)$
\begin{eqnarray}
V\left({\bf d}^3;z\right)=-\sum_{i=1}^3 z^{F({\bf d}^3)+\Sigma_3-a_{ii d_i}}
+z^{J({\bf d}^3)}+\sum_{i=1}^3 z^{a_{ii}d_i}-z^{F({\bf d}^3)+\Sigma_3}\;.
\label{fr3}
\end{eqnarray}
Write equation (\ref{tr19a}) for the 3-dim case
\begin{equation}
\sum_{h\;\in\;\Delta_{{\cal H}}({\bf d}^3)}z^h=\frac{V\left({\bf d}^3;z\right)
}{\left(1-z^{d_1}\right)\left(1-z^{d_2}\right)\left(1-z^{d_3}\right)}\;,
\label{fr4}
\end{equation}
and notice that every term of the generating polynomial at the left hand side
of (\ref{fr4}) has coefficient equal to 1. Moreover, the term of minimal degree
of this polynomial has to coincide with the term having minimal degree in
$V\left({\bf d}^3;z\right)$, and therefore such term has also coefficient
equal to 1. By (\ref{fr3}), polynomial $V\left({\bf d}^3;z\right)$ has four
candidates to have the minimal degree: $z^{J({\bf d}^3)}$ and $z^{a_{ii}d_i}$,
$i=1,2,3$.

We prove that the degrees of the last three terms always exceed $J\left({\bf
d}^3\right)$. Making use of formula (\ref{tr24c}) and relations (\ref{joh1})
and (\ref{eq2}), we obtain the following representation for $J\left({\bf
d}^3\right)$, holding for any permutation $(i,j,k)$ of $(1,2,3)$ :
\begin{eqnarray}
J\left({\bf d}^3\right)=|a_{ij}d_j-a_{ji}d_i|\;.\label{ar5}
\end{eqnarray}
One can always suppose that $a_{ij}>a_{ji}$ up to permutation of indices
$i,j, k$. Since all $a_{ij}$ are positive in the nonsymmetric case, we
have by (\ref{ar5})
\begin{eqnarray}
J\left({\bf d}^3\right)< a_{ij}d_j\;.\label{ar6}
\end{eqnarray}
However, by (\ref{eq2}) we have $a_{ij}d_j=a_{jj}d_j-a_{kj}d_j$, or, in other
words,
\begin{eqnarray}
a_{ij}d_j < a_{jj}d_j\;.\label{ar7}
\end{eqnarray}
Combining the last inequality with (\ref{ar6}) we get, for every $i=1,2,3$,
\begin{eqnarray}
J\left({\bf d}^3\right)< a_{ii}d_i\;.\label{ar8}
\end{eqnarray}
So, $J\left({\bf d}^3\right)$ is the minimal degree of the terms entering
polynomial $V\left({\bf d}^3;z\right)$ and therefore, by comparing the left 
and right hand sides of (\ref{fr4}), we have $\min\Delta_{{\cal H}}\left({\bf 
d}^3\right)=J\left({\bf d}^3\right)$. The second part of (\ref{tr25}) follows 
by (\ref{minmax}).$\;\;\;\;\;\;\Box$

\vspace{.1cm}
As example, consider the semigroup of Figure 1.
\begin{example}$\{d_1,d_2,d_3\}=\{5,11,13\}$\label{ex1}
{\footnotesize \[(a_{ij})=\left(\begin{array}{ccc} 7 & 2 & 1\\ 4 & 3 & 1\\
3 & 1 & 2 \end{array}\right),\;\;\;\left\{\begin{array}{l}a_{11}d_1=
35\\a_{22}d_2=33\\a_{33}d_3=26\end{array}\right.,\;\;\;\left\{\begin{array}{l}
F(5,11,13)=19\\G(5,11,13)=12\\J(5,11,13)=2\end{array}\right.,\;\;\;
\left\{\begin{array}{l}{\sf S}^{\prime}(5,11,13)=\{17,19\}\\
\Delta_{{\cal G}}(5,11,13)=\{1,3,4,6,8,9,14,19\}\\
\Delta_{{\cal H}}(5,11,13)=\{2,7,12,17\}\end{array}\right.,\]
\[H(5,11,13;z)=\frac{Q(5,12,13;z)}{(1-z^{5})
(1-z^{11})(1-z^{13})}\;,\;\;\;\;\;\;\;\;\;\;\;\;\;
Q(5,11,13;z)=1-z^{26}-z^{33}-z^{35}+z^{46}+z^{48}\;,\]
\[\sum_{h\in \Delta_{{\cal H}}} z^h=\frac{ V(5,11,13;z)}{(1-
z^{5})(1-z^{11})(1-z^{13})}\;,\;\;\;\;\;\;\;\;\;V(5,11,13;z)=
z^2-z^{13}-z^{15}-z^{22}+z^{26}+z^{33}+z^{35}-z^{46}.\]}
\end{example}
One more statement can be proven for $\min\Delta_{{\cal H}}({\bf d}^3
)$ if we combine (\ref{tr8}) and (\ref{fr12}) with (\ref{tr24b}) and
(\ref{tr25}):
\begin{eqnarray}
J\left({\bf d}^3\right)\leq c\left({\bf d}^m\right)-G\left({\bf d}^m\right)\;
\;\;\rightarrow\;\;\;J\left({\bf d}^3\right)\leq 1+a_{11}a_{22}a_{33}-\Sigma_3
\;.\label{ar9}
\end{eqnarray}
\section{High -- Dimensional Cases}\label{p4}
In this Section we extend our results about minimal and maximal
values of the $h$--gaps to the numerical semigroups generated by
four and more positive integers. The following Theorem indicates
that for $m>3$ the situation is much more difficult.
\begin{theorem}\label{the5} {\rm (\cite{worm86}, Theorem 3)}\\
The number $\#\{Q\left({\bf d}^m;z\right)\}$ of non--zero coefficients in
the polynomials $Q\left({\bf d}^m;z\right)$ is not bounded by any function
of $m$ for $m\geq 4$, although it is finite for every choice of the generators
$d_i$.
\end{theorem}
However, the next theorem shows how the multiplicity $d_1$ and the embedding
dimension $m$ of a numerical semigroup give an upper bound for
$\#\{Q\left({\bf d}^m;z\right)\}$.
\begin{theorem}\label{the6} {\rm (\cite{fel04}, Theorem 10)}\\
The number of non--zero coefficients in the polynomial $Q\left({\bf d}^m;z
\right),m\geq 4$, is bounded and satisfies the inequality
\begin{eqnarray}
\#\{Q\left({\bf d}^m;z\right)\}\leq d_12^{m-1}-2(m-1)\;.\label{ar10}
\end{eqnarray}
\end{theorem}
Unfortunately, the Frobenius problem for nonsymmetric numerical semigroups
${\sf S}\left({\bf d}^m\right),\;m\geq 4$, is still open nowadays, although
many numerical algorithms are developed \cite{wag05}. However, only the
symmetric semigroups ${\sf S}\left({\bf d}^4\right)$ are to our knowledge
studied \cite{ber79}, \cite{frob94} by means of commutative algebra.

Recall the following basic entity for the $m$-dim semigroups facilitating
further discussion: {\em the Ap\'ery set} of semigroup ${\sf S}\left({\bf
d}^m\right)$ with respect to generator $d_j$ is the set \cite{ros00},
\begin{eqnarray}
Ap\left({\bf d}^m;d_j\right):=\left\{s\in {\sf S}\left({\bf d}^m\right)\quad |
\quad s-d_j\not\in {\sf S}\left({\bf d}^m\right)\right\}\;.\label{ar14}
\end{eqnarray}
Hereafter, we'll make worth of the Ap\'ery set of semigroup with respect to
its multiplicity $d_1$. This set has the following properties
\begin{eqnarray}
\#Ap\left({\bf d}^m;d_1\right)=d_1\;,\;\;\;\;F\left({\bf d}^m\right)=
\max Ap\left({\bf d}^m;d_1\right)-d_1\;,\;\;\;\;0,d_2,d_3\ldots ,d_m\in
Ap\left({\bf d}^m;d_1\right)\;.\label{ar15}
\end{eqnarray}
The generating function ${\cal A}_{d_1}\left({\bf d}^m;z\right)$ for the
Ap\'ery set $Ap\left({\bf d}^m;d_1\right)$ was given in \cite{fel04},
\begin{eqnarray}
{\cal A}_{d_1}\left({\bf d}^m;z\right)=\sum_{s\in Ap({\bf d}^m;d_1)}
z^s=\sum_{k=0}^{d_1-1}z^{k}-(1-z^{d_1})\Phi\left({\bf d}^m;z\right)
\;,\label{ar18}
\end{eqnarray}
where $\Phi\left({\bf d}^m;z\right)$ is the gaps generating function
(\ref{tr3}). Thus, combining (\ref{tr16}) and (\ref{ar18}) we get
\begin{eqnarray}
Q\left({\bf d}^m;z\right)={\prod_{j=2}^m(1-z^{d_j})}\cdot {\cal A}_
{d_1}\left({\bf d}^m;z\right)\;.\label{ar21}
\end{eqnarray}
Notice that (\ref{ar21}) and the last containment in (\ref{ar15}) are
sufficient to prove Theorem \ref{the6}.

Denote by $D_i\left({\bf d}^m;d_1\right)$, $1\leq i\leq m-1$, the
sets of integers
\begin{eqnarray}
D_i({\bf d}^m;d_1)&=&\left\{d_{k_1}+\ldots +d_{k_i}\;|\;k_1,\ldots k_i
\in\{2,\ldots m\}\;,\;\;\;k_1\neq k_2\neq \ldots\neq k_i\right\}\;,\;
\mbox{i.e.}\label{ar210}\\
D_1({\bf d}^m;d_1)&=&\left\{d_k\;|\;k\in\{2,\ldots m\}\right\}\;,
\nonumber\\
D_2({\bf d}^m;d_1)&=&\left\{d_k+d_l\;|\;k,l\in\{2,\ldots m\},\;
k\neq l\right\}\;,\;\ldots \;,\;\;\mbox{etc}.\nonumber
\end{eqnarray}
Finally, the last set $D_{m-1}\left({\bf d}^m;d_1\right)=\left\{\Sigma_m-d_1
\right\}$ consists of only one integer.

Construct, for every $i=1,\ldots, m-1$, the direct sum of the sets
\begin{eqnarray}
Ap\left({\bf d}^m;d_1\right)\oplus D_i\left({\bf d}^m;d_1\right)=\left\{\alpha+
\delta_i\;|\;\alpha\in Ap\left({\bf d}^m;d_1\right),\;\delta_i\in
D_i\left({\bf d}^m;d_1\right)\right\}\;,\label{ar211}
\end{eqnarray}
and denote by $B_i\left({\bf d}^m\right)$ the set of degrees of the terms
$z^{C_{j,i}}$ for the syzygies of the $i$th kind in (\ref{bet1}). By comparison
of (\ref{ar21}) and (\ref{bet1}) we have
\begin{eqnarray}
B_i\left({\bf d}^m\right)\subset Ap\left({\bf d}^m;d_1\right)\oplus
D_i\left({\bf d}^m;d_1\right)\;,\label{ar212}
\end{eqnarray}
that gives immediately the upper bounds for the Betti numbers, $\beta_i\left(
{\bf d}^m\right)\leq d_1{m-1\choose i}$, which satisfy (\ref{ar10}). In
particular,
\begin{eqnarray}
B_{m-1}\left({\bf d}^m\right)\subset Ap\left({\bf d}^m;d_1\right)\oplus
\left\{\Sigma_m-d_1\right\}\;.\label{ar213}
\end{eqnarray}
Notice that the elements of the set $B_{m-1}\left({\bf d}^m\right)$ are those
elements of the set $Ap\left({\bf d}^m;d_1\right)\oplus\left\{\Sigma_m-d_1
\right\}$ which are survived in $Q\left({\bf d}^m;z\right)$ after all
cancellations in (\ref{ar21}).

In order to prove the  main result of this Section (next Theorem \ref{the7}) we
investigate the relationship between two sets, $B_{m-1}\left({\bf d}^m\right)$
and ${\sf S}^{\prime}\left({\bf d}^m\right)$, which have equal cardinalities
(see \ref{bet2n}).

Enumerate the elements $F_i$ of the set ${\sf S}^{\prime}\left({\bf d}^m
\right)$ in ascending order, $F_{i-1}<F_i$, $i=1,\ldots, t\left({\bf d}^m
\right)$, and write, according to (\ref{tr6a}),
\begin{eqnarray}
F_j+d_1\in {\sf S}\left({\bf d}^m\right)\;,\;\;\;F_j\in\Delta\left({\bf
d}^m\right)\;.\label{ar21a}
\end{eqnarray}
Simple comparison of (\ref{ar21a}) with (\ref{ar14}) gives
\begin{eqnarray}
F_j+d_1\in Ap({\bf d}^m;d_1)\;,\;\;\;j=1,\ldots, t({\bf d}^m
)\;,\label{ar21b}
\end{eqnarray}
that implies
\begin{eqnarray}
{\sf S}^{\prime}\left({\bf d}^m\right)\oplus\left\{\Sigma_m\right\}\subset
Ap\left({\bf d}^m;d_1\right)\oplus\left\{\Sigma_m-d_1\right\}\;.\label{ar214}
\end{eqnarray}
Consider what kind of elements of the set ${\sf S}^{\prime}\left({\bf d}^m
\right)\oplus\left\{\Sigma_m\right\}$ is survived in $Q\left({\bf d}^m;z\right)$
after all cancellations in (\ref{ar21}).
\begin{lemma}\label{lem1}{\rm $\;$}\\
The set of values $F_j+\Sigma_m$ coincides with the set $B_{m-1}\left({\bf
d}^m\right)$,
\begin{eqnarray}
{\sf S}^{\prime}\left({\bf d}^m\right)\oplus\left\{\Sigma_m\right\}=
B_{m-1}\left({\bf d}^m\right)\;\;.\label{ar216}
\end{eqnarray}
\end{lemma}
{\sf Proof}$\;\;\;$First, prove that elements $F_j+\Sigma_m$ are
not cancelled in (\ref{ar21}), or, in other words, \\ ${\sf
S}^{\prime}\left({\bf d}^m\right)\oplus\left\{\Sigma_m\right\}\subseteq
B_{m-1}\left({\bf d}^m\right)$. According to (\ref{ar21}) and (\ref{ar21b}),
an element $F_j+\Sigma_m$ is the degree of a term of polynomial $Q\left({\bf
d^m};z\right)$ iff it is not cancelled by a term (with opposite sign) of
degree $y+\Sigma_m-\Sigma_ {odd}$, where $y+d_1\in Ap\left({\bf d}^m;d_1
\right)$ and  $\Sigma_{odd}$ is the sum of an odd number of generators
different from $d_1$. Note that, by the same definition (\ref{ar14}), $y\in
\Delta\left({\bf d}^m\right)$. Let, by way of contradiction, $F_j+\Sigma_m$ be
cancelled for some $j$. Then there exists $y$ such that
\begin{eqnarray}
F_j+\Sigma_m=y+\Sigma_m-\Sigma_{odd}\;\;\;\;\longrightarrow\;\;\;\;
y=F_j+\Sigma_{odd}\;.\label{ar21c}
\end{eqnarray}
However, the last equality in (\ref{ar21c}) contradicts our assumptions:
$y\in \Delta\left({\bf d}^m\right)$ and $F_j\in {\sf S}^{\prime}\left({\bf
d}^m\right)$. Thus we have
\begin{eqnarray}
{\sf S}^{\prime}\left({\bf d}^m\right)\oplus\left\{\Sigma_m\right\}\subseteq
B_{m-1}\left({\bf d}^m\right)\;.\label{ar215}
\end{eqnarray}
Finally, combining equality (\ref{bet2n}) of both cardinalities, $\#{\sf S}^{
\prime}\left({\bf d}^m\right)$ and $\#B_{m-1}\left({\bf d}^m\right)$, allows
us to conclude the proof of Lemma \ref{lem1}, so obtaining the bijection :
\begin{eqnarray}
C_{j,m-1}=F_j+\Sigma_m\;,\;\;\;j=1,\ldots, t\left({\bf d}^m\right)\;.
\label{ar217}
\end{eqnarray}
$\;\;\;\;\;\;\Box$
\vspace{.1cm}

\begin{lemma}\label{lem02}{\rm $\;$}\\
The elements of the set $B_{m-1}\left({\bf d}^m\right)$ do not coincide with 
elements of any other set $B_{m-i}\left({\bf d}^m\right)$, $i\geq 2$.
\end{lemma}
{\sf Proof}$\;\;\;$  By Lemma \ref{lem1}, the elements of the set
$B_{m-1}\left( {\bf d}^m\right)$ are given by (\ref{ar216}).
Because of cancellations  in (\ref{ar21}), an element of the set
$B_{m-1}\left({\bf d}^m\right)$ could coincide only with an
element of a set $B_{m-i}\left({\bf d}^m\right)$, for $i$ odd,
$i\geq 3$.  To exclude this,  it is thus sufficient to repeat the
proof, by way of contradiction, of the first part of Lemma
\ref{lem1}, substituting $\Sigma_{odd}$ by $\Sigma_{even}$.
$\;\;\;\;\;\;\Box$

Below we give a generalization of Theorem \ref{the4} for $m\geq 4$. Denote by
$\Xi\left({\bf d}^m\right)$ the union of the sets $B_i\left({\bf d}^m\right)$
excluding its maximal element $\deg Q\left({\bf d}^m;z\right)$, i.e.
\begin{eqnarray}
\Xi\left({\bf d}^m\right)=\cup_{i=1}^{m-1}B_i\left({\bf d}^m\right)\setminus\{
\deg Q\left({\bf d}^m;z\right)\}\;,\label{ar218}
\end{eqnarray}
where $\#\Xi\left({\bf d}^m\right)=\#\{Q\left({\bf d}^m;z\right)\}-2$.
Define also the complementary set $\overline\Xi\left({\bf d}^m\right)$,
\begin{eqnarray}
\overline\Xi\left({\bf d}^m\right):=\{\deg Q\left({\bf d}^m;z\right)-\xi_k,
\;| \;\xi_k\in \Xi\left({\bf d}^m\right)\}\;,\;\;\;\;\#\Xi\left({\bf d}^m
\right)=\#\overline\Xi\left({\bf d}^m\right)\;.\label{ar225}
\end{eqnarray}
\begin{lemma}\label{lemv}{\rm $\;$}\\
The set  $\Theta\left({\bf d}^m\right)$ of the degrees of all terms entering
polynomial $V\left({\bf d}^m;z\right)$ satisfies:
\begin{eqnarray}
\Theta\left({\bf d}^m\right)\subseteq\Xi\left({\bf d}^m\right)\cup\overline
\Xi\left({\bf d}^m\right)\;.\label{eqv}
\end{eqnarray}
\end{lemma}
{\sf Proof}$\;\;\;$Combining (\ref{bet1}) and (\ref{tr18}) we obtain
\begin{eqnarray}
V\left({\bf d}^m;z\right)=\sum_{k=1}^{m-1}(-1)^{m-k+1}\sum_{j=1}^{\beta_k
\left({\bf d}^m\right)}\left\{z^{\deg Q\left({\bf d}^m;z\right)-C_{j,k}}+(-1)^
mz^{C_{j,k}}\right\}-\label{tr19}\\
\left\{1+(-1)^mz^{\deg Q({\bf d}^m;z)}\right\}.\nonumber
\end{eqnarray}
The two last terms are cancelled with the similar terms existing in the
left sums for $k=m-1,\;j=\beta_{m-1}\left({\bf d}^m\right)$. The degrees of the
remained terms in (\ref{tr19}) are coming from both sets, $\Xi\left({\bf d}^m
\right)$ and $\overline\Xi\left({\bf d}^m\right)$, that proves Lemma.
$\;\;\;\;\;\;\Box$
\begin{theorem}\label{the7}{\rm $\;$}\\
Let a nonsymmetric semigroup ${\sf S}\left({\bf d}^m\right)$ be given, and
the Hilbert series $H\left({\bf d}^m;z\right)$ of its graded subring be
given by (\ref{tr16}), (\ref{bet1}). Then $\min\Delta_{{\cal H}}\left({\bf
d}^m\right)$ is equal to the difference of the degrees in the two last terms
of polynomial $Q\left({\bf d}^m;z\right)$,
\begin{eqnarray}
\min\Delta_{{\cal H}}\left({\bf d}^m\right)=\deg Q\left({\bf d}^m;z\right)-
\max\Xi\left({\bf d}^m\right)\;.\label{ar12}
\end{eqnarray}
\end{theorem}
To prove this theorem we need the auxiliary Lemma \ref{lem2}.
\begin{lemma}\label{lem2}{\rm $\;$}\\
Let a nonsymmetric semigroup ${\sf S}\left({\bf d}^m\right)$ be given. Then
\begin{eqnarray}
F_{t({\bf d}^m)}=F({\bf d}^m),\;\;\;{\rm and}\;\;\;F_{t({\bf d}^m)-1}=
\max\Delta_{{\cal H}}\left({\bf d}^m\right)\;.\label{ar22}
\end{eqnarray}
\end{lemma}
{\sf Proof} $\;\;\;$By definition (\ref{tr4}), the Frobenius number $F\left({
\bf d}^m\right)$ is the maximal value in $\Delta\left({\bf d}^m\right)$, and
therefore, in accordance with (\ref{tr6a}), $F\left({\bf d}^m\right)=\max 
{\sf S}^{\prime}\left({\bf d}^m\right)$. Hence, $F_{t({\bf d}^m)}=
F\left({\bf d}^m\right)$.

Let the element $F_{t({\bf d}^m)-1}\in {\sf S}^{\prime}\left({\bf d}^m\right)$ 
be given. Consider the difference $x=F\left({\bf d}^m\right)-F_{t({\bf d}^m)-
1}$, and prove that $x\in\Delta\left({\bf d}^m\right)$. Suppose, by way of 
contradiction, that $x\in {\sf S}\left({\bf d}^m\right)$. However, we have 
$F_{t({\bf d}^m)-1}+x=F\left({\bf d}^m\right)\in\Delta\left({\bf d}^m\right)$, 
that contradicts our assumption, $F_{t({\bf d}^m)-1}\in {\sf S}^{\prime}
\left({\bf d}^m\right)$. Hence $F\left({\bf d}^m\right)-F_{t({\bf d}^m)-1}\in
\Delta\left({\bf d}^m\right)$, that, according to (\ref{tr9}), yields
\begin{eqnarray}
F_{t({\bf d}^m)-1}\in\Delta_{{\cal H}}\left({\bf d}^m\right)\;.\label{ar221}
\end{eqnarray}

Prove that $F_{t({\bf d}^m)-1}$ is the maximal value in $\Delta_{{\cal H}}
\left({\bf d}^m\right)$. Suppose, by way of contradiction, that there exists a
value $y_1\in\Delta_{{\cal H}}\left({\bf d}^m\right)$ such that $F_{t({\bf d}^
m)-1}<y_1<F_{t({\bf d}^m)}$. Then, there exists necessarily $y_2\in {\sf S}
\left({\bf d}^m\right)$ such that $y_1+y_2\in\Delta\left({\bf d}^m\right)$. 
Summarize the known facts as follows :
\begin{eqnarray}
&&y_1\in\Delta_{{\cal H}}\left({\bf d}^m\right)\;,\;\;\;F\left({\bf d}^m
\right)-y_1\in \Delta_{{\cal H}}\left({\bf d}^m\right)\;,\label{ar222}\\
&&y_2\in {\sf S}\left({\bf d}^m\right)\;,\;\;\;y_1+y_2\in\Delta\left({\bf
d}^m\right)\;.\label{ar223}
\end{eqnarray}
First, suppose that $y_1+y_2\not\in\Delta_{{\cal H}}\left({\bf d}^m\right)$, 
that implies
\begin{eqnarray}
F\left({\bf d}^m\right)-(y_1+y_2)=(F({\bf d}^m)-y_1)-y_2=s\in
{\sf S}\left({\bf d}^m\right)\;\;\rightarrow\;\;F({\bf d}^m)-y_1=
y_2+s\in {\sf S}\left({\bf d}^m\right)\nonumber
\end{eqnarray}
and contradicts (\ref{ar222}).

Next, suppose that $y_1+y_2\in\Delta_{{\cal H}}\left({\bf d}^m\right)$.
However, this means that we have found {\em another value} $y_1+y_2$ which
satisfies
\begin{eqnarray}
F_{t({\bf d}^m)-1}<y_1+y_2<F_{t({\bf d}^m)}\;.\nonumber
\end{eqnarray}
Further repetitions of this procedure contradict the fact that the interval
$\left[F_{t({\bf d}^m)-1},F_{t({\bf d}^m)}\right]$ is finite and makes our
assumption about a value $y_1\in\Delta_{{\cal H}}\left({\bf d}^m\right)$
such that $F_{t({\bf d}^m)-1}<y_1<F_{t({\bf d}^m)}$ wrong. Thus, the 2nd
equality in (\ref{ar22}) is proven.
$\;\;\;\;\;\;\Box$

\vspace{.1cm}
To complete the proof of Theorem \ref{the7}, observe that, combining the 1st
equality in (\ref{ar22}) with equality (\ref{ar217}), for $j=t\left({\bf d}^m
\right)$, we obtain:
\begin{equation}
F\left({\bf d}^m\right)=\deg Q\left({\bf d}^m;z\right)-\Sigma_m\;,\label{FQS}
\end{equation}
(being $\deg Q\left({\bf d}^m;z\right)=C_{\beta_{m-1},m-1}$).

Using equalities (\ref{minmax}) and (\ref{FQS}) we get
\begin{equation}
\min\Delta_{{\cal H}}\left({\bf d}^m\right)=\deg Q\left({\bf d}^m;z\right)-
\Sigma_m-\max\Delta_{{\cal H}}\left({\bf d}^m\right)\;.\label{mindelta}
\end{equation}
Combining  the 2nd equality in (\ref{ar22}) with (\ref{ar217}), we get
\begin{equation}
\max\Delta_{{\cal H}}\left({\bf d}^m\right)=\max\left\{B_{m-1}\left({\bf
d}^m\right)\setminus\{\deg Q\left({\bf d}^m;z\right)\}\right\}-\Sigma_m\;.
\label{maxdelta}
\end{equation}
On the other hand, by Theorem \ref{the3} and Lemma \ref{lemv},
\begin{eqnarray}
\max\Delta_{{\cal H}}\left({\bf d}^m\right)+\Sigma_m\;\;\;\mbox{is either
equal to}\;\;\;\max\Xi\left({\bf d}^m\right)\;\;\mbox{or to}\;\;\max\overline
\Xi\left({\bf d}^m\right)\;.\nonumber
\end{eqnarray}
Equation (\ref{maxdelta}) says that $\max\Delta_{{\cal H}}\left({\bf d}^m
\right)+\Sigma_m$ belongs to set $\Xi\left({\bf d}^m\right)$, and, more
precisely:
\begin{equation}
\max\Delta_{{\cal H}}\left({\bf d}^m\right)= \max\Xi\left({\bf
d}^m\right)-\Sigma_m\;.\label{maxxi}
\end{equation}
Combining (\ref{mindelta}) with (\ref{maxxi}) we obtain the statement of the
theorem. $\;\;\;\;\;\;\Box$

As a spin-off of the proof of Theorem \ref{the7} we have a
corollary concerning the maximal degrees of the terms entering the
syzygies of different kinds in (\ref{bet1}).
\begin{corollary}\label{cor02}{\rm $\;$}\\
\begin{equation}
\max\left\{B_{m-1}\left({\bf d}^m\right)\setminus\{\deg Q\left({\bf d}^m;z
\right)\}\right\}>\max\left\{B_{m-i}\left({\bf d}^m\right)\right\}\;,\;\;\;
i\geq 2\;.\label{maxsy1}
\end{equation}
\end{corollary}
{\sf Proof} $\;\;\;$By consequence of (\ref{maxdelta}) and (\ref{maxxi}),
\begin{equation}
\max\Xi\left({\bf d}^m\right)=\max\left\{B_{m-1}\left({\bf d}^m\right)
\setminus\{\deg Q\left({\bf d}^m;z\right)\}\right\}\;.\label{maxsy2}
\end{equation}
However by definition (\ref{ar218}) of the set $\Xi\left({\bf d}^m\right)$, 
we have also $\max\Xi\left({\bf d}^m\right)\geq\max\left\{B_{m-i}\left({\bf 
d}^m\right)\right\}$, $i\geq 2$. On the other hand, by consequence of Lemma 
\ref{lem02}, the last inequality must be strict since $\max\left\{ B_{m-1}
\left({\bf d}^m\right)\setminus\{\deg Q\left({\bf d}^m;z\right)\}\right\}\not
\in B_{m-i}\left({\bf d}^m\right)$ for any $i\geq 2$. Combining both these 
facts, we come to (\ref{maxsy1}).$\;\;\;\;\;\;\Box$
\begin{lemma}\label{lem3}{\rm $\;$}\\
Let a nonsymmetric semigroup ${\sf S}\left({\bf d}^m\right)$ be given. Then
\begin{eqnarray}
{\sf S}^{\prime}\left({\bf d}^m\right)\setminus\left\{F\left({\bf d}^m\right)
\right\}\subseteq \Delta_{{\cal H}}\left({\bf d}^m\right)\;.\label{ar229}
\end{eqnarray}
\end{lemma}
{\sf Proof} $\;\;\;$ Let  $F_j\in {\sf S}^{\prime}\left({\bf d}^m\right)$,
$j\not=t\left({\bf d}^m\right)$. By definition (\ref{tr6a}), it belongs to
the set $\Delta\left({\bf d}^m\right)$. Suppose, by way of contradiction,
that $F_j\not\in\Delta_{{\cal H}}\left({\bf d}^m\right)$. Then, there exists
an element $s=F\left({\bf d}^m\right)-F_j$ such that $s\in {\sf S}\left(
{\bf d}^m\right)$, and therefore $F_j+s=F\left({\bf d}^m\right)\in\Delta
\left({\bf d}^m\right)$. This contradicts the assumption that $F_j$ belongs 
to $S^{\prime}\left({\bf d}^m\right)$. Hence $F_j\in\Delta_{{\cal H}}\left(
{\bf d}^m\right)$.$\;\;\;\;\;\;\Box$

Lemma \ref{lem3} has an important Corollary which would be made worth in
Section \ref{p5}.
\begin{corollary}\label{cor2a}{\rm $\;$}\\
\begin{eqnarray}
\max\Delta_{{\cal H}}\left({\bf d}^m\right)-\min\Delta_{{\cal H}}\left({\bf
d}^m\right)\geq t\left({\bf d}^m\right)-2\;,\label{ar228}
\end{eqnarray}
\begin{eqnarray}
\max\Delta_{{\cal H}}\left({\bf d}^m\right)\geq \frac{1}{2}\left[F\left(
{\bf d}^m\right)+t\left({\bf d}^m\right)\right]-1\;,\;\;\;
\min\Delta_{{\cal H}}\left({\bf d}^m\right)\leq\frac{1}{2}\left[F\left(
{\bf d}^m\right)-t\left({\bf d}^m\right)\right]+1\;.\;\;\label{ar228a}
\end{eqnarray}
\end{corollary}
{\sf Proof} $\;\;\;$The proof is based on Lemma \ref{lem3} and follows from
a chain of inequalities. Taking in mind $\#\left\{{\sf S}^{\prime}\left({\bf
d}^m\right)\setminus \left\{F\left({\bf d}^m\right)\right\}\right\}=t\left(
{\bf d}^m\right)-1$, we have
\begin{eqnarray}
\max\Delta_{{\cal H}}\left({\bf d}^m\right)-\min\Delta_{{\cal H}}\left({\bf
d}^m\right)+1\geq \#\Delta_{{\cal H}}\left({\bf d}^m\right)\geq t\left({\bf
d}^m\right)-1\;,\nonumber
\end{eqnarray}
that proves (\ref{ar228}). Combining now  (\ref{minmax}) with (\ref{ar228}) 
we arrive at (\ref{ar228a})$\;\;\;\;\;\;\Box$

One more Corollary deals with the maximal and minimal values of the set
$\Xi\left({\bf d}^m\right)$.
\begin{corollary}\label{cor2}{\rm $\;$}\\
\begin{eqnarray}
\min\Xi\left({\bf d}^m\right)+\max\Xi\left({\bf d}^m\right)>\deg Q\left({\bf
d}^m;z\right)\;.\label{ar224}
\end{eqnarray}
\end{corollary}
{\sf Proof} $\;\;\;$ By Theorem \ref{the3} the minimal degree among the terms 
which contribute to $V\left({\bf d}^m;z\right)$ is equal to $\min\Delta_{{\cal 
H}}\left({\bf d}^m\right)$. On the other hand, by Theorem \ref{the7} this 
degree is given by $\deg Q\left({\bf d}^m;z\right)-\max\Xi\left({\bf d}^m
\right)$. Since the degree $\min\Xi\left({\bf d}^m\right)$ contributes 
also to the polynomial $V\left({\bf d}^m;z\right)$, we conclude
\begin{eqnarray}
\deg Q\left({\bf d}^m;z\right)-\max\Xi\left({\bf d}^m\right)<\min\Xi\left(
{\bf d}^m\right)\;,\nonumber
\end{eqnarray}
that proves Corollary.$\;\;\;\;\;\;\Box$

\vspace{.2cm}
We finish this Section with one Example for the 4-dim numerical semigroup
studied in \cite{kill00}. The polynomial $Q\left({\bf d}^m;z\right)$ was
obtained by means of diagrammatic calculation on the set $\Delta\left({\bf 
d}^m\right)$ developed for $m=3$ and extended for $m\geq 4$ \cite{fel04}.
\begin{example}$\{d_1,d_2,d_3,d_4\}=\{103,133,165,228\}$\label{ex2}
{\footnotesize
\begin{eqnarray}
\left\{\begin{array}{lll}
a_{11}d_1&=&824\\a_{22}d_2&=&1197\\a_{33}d_3&=&825\\a_{44}d_4&=&1368
\end{array}\right.,\;\;\;\;\;\left\{\begin{array}{lrl}F(103,133,165,228)&=&
1436\\G(103,133,165,228)&=&840\\
\Delta_{{\cal H}}(103,133,165,228)&=&\{19,\ldots,1417\}
\\ \#\left\{Q(103,133,165,228;z)\right\}&=&26\end{array}\right.,\;\;\;\;\;
\left.\begin{array}{lrl}\deg Q&=&2065\\\max\Xi&=&2046\\\min\Xi&=&824\\
\min\Delta_{{\cal H}}&=&2065-2046\end{array}\right.\nonumber
\end{eqnarray}
\begin{eqnarray}
{\sf S}^{\prime}(103,133,165,228)&=&\{1145,1316,1355,1374,1417,1436
\}\;,\;\;\;\;\beta_1=7\;,\;\;\;\beta_2=12\;,\;\;\;\beta_3=6\;,\nonumber\\
H(103,133,165,228;z)&=&\frac{Q(103,133,165,228;z)}{(1-
z^{103})(1-z^{133})(1-z^{165})(1-z^{228})}\;,\;\;\;\;\Sigma_4=629\;,\nonumber\\
Q(103,133,165,228;z)&=&1-z^{824}-z^{825}-z^{1077}-z^{1096}-
z^{1197}-z^{1216}+z^{1319}+z^{1362}-z^{1368}+z^{1489}+\nonumber\\
&&z^{1508}+z^{1533}+z^{1546}+z^{1609}+z^{1737}+z^{1756}-z^{1774}+z^{1780}+
z^{1881}+z^{1900}-\nonumber\\
&&z^{1945}-z^{1984}-z^{2003}-z^{2046}-z^{2065}\;\nonumber
\end{eqnarray}}\end{example}
\section{Special Kinds of Numerical Semigroups}\label{p5}
In this Section we focus on three special kinds of numerical semigroups, namely,
the pseudo--symmetric semigroups, semigroups with maximal embedding dimension,
and semigroups of maximal and almost maximal length.
\subsection{Pseudo--symmetric semigroups}\label{s1p5}
Following \cite{baru97} and \cite{frob87} call the nonsymmetric semigroup
${\sf S}\left({\bf d}^m\right)$ {\em pseudo--symmetric} if
\begin{eqnarray}
{\sf S}^{\prime}\left({\bf d}^m\right)=\left\{\frac{F\left({\bf d}^m\right)}{2}
,F\left({\bf d}^m\right)\right\}\;.\label{fr12a}
\end{eqnarray}
Notice that $F\left({\bf d}^m\right)$ is necessarily an even number and
$t\left({\bf d}^m\right)=2$. However the last equality alone is not enough to 
provide the pseudo--symmetricity of semigroup ${\sf S}\left({\bf d}^m\right)$.
\begin{corollary}\label{cor3}{\rm $\;$}\\
Every pseudo--symmetric semigroup ${\sf S}\left({\bf d}^m\right)$ satisfies:
\begin{eqnarray}
\Delta_{{\cal H}}\left({\bf d}^m\right)=\left\{\frac{F\left({\bf d}^m\right)}
{2}\right\}\;.\label{ar31}
\end{eqnarray}
\end{corollary}
{\sf Proof} $\;\;\;$By Lemma \ref{lem2} we have $\max\Delta_{{\cal H}}\left(
{\bf d}^m\right)=\frac{1}{2}F\left({\bf d}^m\right)$. On the other hand,
according (\ref{minmax}), we get $\min\Delta_{{\cal H}}\left({\bf d}^m\right)
=\frac{1}{2}F\left({\bf d}^m\right)$. Thus, $\#\Delta_{{\cal H}}\left({\bf
d}^m\right)=1$ and (\ref{ar31}) is proven.$\;\;\;\;\;\;\Box$
\begin{remark}\label{re1}
{\rm Notice that for all $m>3$ there exist semigroups with $\#\Delta_{{\cal H}}
\left({\bf d}^m\right)=1$ which are not pseudo--symmetric. However, the case 
$m=3$ is a special one.}
\end{remark}
\begin{corollary}\label{cor4}{\rm $\;$}\\
Let a nonsymmetric semigroup ${\sf S}\left({\bf d}^3\right)$ be given, and
$\#\Delta_{{\cal H}}\left({\bf d}^3\right)=1$. Then ${\sf S}\left({\bf
d}^3\right)$ is pseudo--symmetric.
\end{corollary}
{\sf Proof} $\;\;\;$ By definition of the set $\Delta_{{\cal H}}\left({\bf 
d}^m\right)$, the equality $\#\Delta_{{\cal H}}\left({\bf d}^m\right)=1$ 
implies that $\min\Delta_{{\cal H}}\left({\bf d}^m\right)=\max\Delta_{{\cal H}}
\left({\bf d}^m\right)=\frac{1}{2}F\left({\bf d}^m\right)$. By Lemma 
\ref{lem2}, this value is equal to $F_{t({\bf d}^m)-1}$. Taking in mind that 
for the 3-dim nonsymmetric semigroups $t\left({\bf d}^3\right)=2$, we conclude 
that ${\sf S}\left({\bf d}^3\right)$ is pseudo--symmetric.$\;\;\;\;\;\;\Box$

\vspace{.2cm}
In the rest of this Section we study the pseudo--symmetric semigroups ${\sf S}
\left({\bf d}^3\right)$ in more details. We prove that all semigroups ${\sf S}
\left({\bf d}^3\right)$ with one $h$--gap are generated by a 3--parametric
family of triples $\{d_1,d_2,d_3\}$ and give their parametric representation.

Recall the following basic representations of the conductor $c\left({\bf d}^3
\right)$, the minimal value in $\Delta_{{\cal H}}\left({\bf d}^3\right)$, the 
cardinality $\#\Delta_{{\cal H}}\left({\bf d}^3\right)$ and generators 
$d_1,d_2$ and $d_3$ through the matrix elements $a_{ij}$ \cite{fel04} :
\begin{eqnarray}
c\left({\bf d}^3\right)&=&1+(a_{21}+a_{31})(a_{12}+a_{32}(a_{31}+a_{32})-
(a_{12}a_{23}+a_{23}a_{31}+a_{31}a_{12})-\label{dub1}\\
&&(a_{21}a_{13}+a_{13}a_{32}+a_{32}a_{21})-(a_{31}a_{32}+a_{12}a_{13}+a_{32}
a_{21})+\nonumber\\
&&max\{a_{12}a_{23}a_{31},a_{21}a_{13}a_{32}\}\;,\nonumber
\end{eqnarray}
\begin{eqnarray}
\min\Delta_{{\cal H}}\left({\bf d}^3\right)=|a_{12}a_{23}a_{31}-a_{21}a_{13}
a_{32}|\;,\;\;\;\;\;\;\;d_1=a_{12}a_{23}+a_{32}a_{13}+a_{12}a_{13}\;,
\label{dub2}\\
\#\Delta_{{\cal H}}\left({\bf d}^3\right)=\min\left\{a_{12}a_{23}a_{31},
a_{21}a_{13}a_{32}\right\}\;,\;\;\;\;\;\;\;d_2=a_{23}a_{31}+a_{13}a_{21}+
a_{23}a_{21}\;,\label{dub3}\\
d_3=a_{31}a_{12}+a_{21}a_{32}+a_{31}a_{32}\;.\label{dub4}
\end{eqnarray}
By definition of pseudo--symmetric semigroups and the 1st formula in
(\ref{dub3}) the elements of the Johnson matrix $(a_{ij})$ satisfy
\begin{eqnarray}
a_{12}=a_{23}=a_{31}=1\;,\;\;a_{21}a_{13}a_{32}\geq 2\;,\label{dub5}
\end{eqnarray}
since at least one of three matrix elements $a_{21},a_{13},a_{32}$ exceeds 1.
According to the 2nd formulas in (\ref{dub2}), (\ref{dub3}) and (\ref{dub4})
the generating triple ${\bf d}^3$ is given by :
\begin{eqnarray}
d_1=1+a_{32}a_{13}+a_{13}\;,\;\;\;\;\;d_2=1+a_{13}a_{21}+a_{21}\;,\;\;\;
\;\;d_3=1+a_{21}a_{32}+a_{32}\;.\label{dub6}
\end{eqnarray}
The straightforward calculation of $c\left({\bf d}^3\right)$ and $G\left({\bf 
d}^3\right)$ for pseudo--symmetric semigroup ${\sf S}\left({\bf d}^3\right)$ 
gives
\begin{eqnarray}
&&c\left({\bf d}^3\right)=2a_{21}a_{13}a_{32}-1\;,\;\;\;\;\;
G\left({\bf d}^3\right)=a_{21}a_{13}a_{32}\;,\;\;\;\;\;\#\Delta_{{\cal G}}
\left({\bf d}^3\right)=a_{21}a_{13}a_{32}-1\;,\nonumber\\
&&\min\Delta_{{\cal H}}\left({\bf d}^3\right)=\max\Delta_{{\cal H}}\left({\bf 
d}^3\right)=a_{21}a_{13}a_{32}-1\;.\label{dub7}
\end{eqnarray}
Note that the formula for $\min\Delta_{{\cal H}}\left({\bf d}^3\right)$ in
(\ref{dub2}) is in full agreement with (\ref{fr12a}) and (\ref{dub7}).

Equations (\ref{dub6}) define a nondegenerate mapping $\{a_{32},a_{13},a_{21}\}
\longmapsto\{d_1,d_2,d_3\}$ since its Jacobian does not vanish,
\begin{eqnarray}
\det\left(\frac{\partial^2 d_i}{\partial a_{jk}\;\partial a_{pq}}\right)=1+
a_{32}+a_{13}+a_{21}+a_{32}a_{13}+a_{13}a_{21}+a_{21}a_{32}+2a_{32}a_{13}a_{21}
\geq 14\;.\nonumber
\end{eqnarray}
Therefore all different (up to permutation) triples $\{a_{32},a_{13},a_{21}\}$ 
give rise to different (up to permutation) triples $\{d_1,d_2,d_3\}$.

Summarizing (\ref{dub5}), (\ref{dub6}) and Corollary \ref{cor4} we formulate 
the following Theorem.
\begin{theorem}\label{the7a}{\rm $\;$}\\
A nonsymmetric semigroup ${\sf S}\left({\bf d}^3\right)$ is pseudo--symmetric 
iff there exist three positive integers $a$, $b$ and $c$ ( $a b c\geq 2$ ) 
such that ${\bf d}^3$ is given by:
\begin{eqnarray}
d_1=1+ab+b\;,\;\;\;\;\;d_2=1+bc+c\;,\;\;\;\;\;d_3=1+ca+a\;.\label{ar227}
\end{eqnarray}
\end{theorem}
\subsection{Semigroups with Maximal Embedding Dimension}\label{s2p5}
In this Section we apply our results of Section \ref{p4} to numerical
semigroups of a special kind, ${\sf S}\left({\bf d}^m_{MED}\right)$ where
${\bf d}^m_{MED}=(m,d_2\ldots,d_m)$, called \cite{sal79} semigroups with {\em
maximal embedding dimension} (MED). Their study is motivated by the following 
Theorem.
\begin{theorem}{\rm (\cite{heku71}, Proposition 1.14)}\label{the8}\\
The numerical semigroup ${\sf S}\left({\bf d}^m_{MED}\right)$, $m\geq 3$,
is never symmetric.
\end{theorem}
There are many known results on the MED--semigroups \cite{heku71}, 
\cite{sal79} :
\begin{eqnarray}
Ap\left({\bf d}^m_{MED};m\right)=\{0,d_2,\ldots ,d_m\}\;,\;\;\;\#Ap\left({\bf
d}^m_{MED};m\right)=m\;,\;\;\;t\left({\bf d}^m_{MED}\right)=m-1\;,\;\;
\label{ar24y}\\
F\left({\bf d}^m_{MED}\right)=d_m-m\;,\;\;\;G\left({\bf d}^m_{MED}\right)=
\frac{1}{m}\sum_{k=2}^md_k-\frac{m-1}{2}\;,\;\;\;\beta_k\left({\bf d}^m_{MED}
\right)=k{m\choose k+1}\;.\label{ar24}
\end{eqnarray}
Note that the divisibility by $m$ of the sum $\sum_{k=2}^md_k$ is guaranteed
by the minimality of the set of generators, since $d_i\neq d_j\pmod m$,
$i\neq j\in\{2,\ldots,m\}$. The minimality of the set of generators implies
also that
\begin{eqnarray}
d_m\leq (m-1)d_2-m\;.\label{ar24a}
\end{eqnarray}
The sum of all Betti numbers of a MED--semigroup is less than the upper
bound in Theorem \ref{the6}
\begin{eqnarray}
\sum_{k=0}^{m-1}\beta_k\left({\bf d}^m_{MED}\right)=(m-2)2^{m-1}+2<m\;2^{m-1}-
2(m-1)\;.\nonumber
\end{eqnarray}
The set ${\sf S}^{\prime}\left({\bf d}^m_{MED}\right)$ is given by $\left\{
d_2-m,\ldots,d_{m-1}-m,d_m-m\right\}$.

Simple structure of the MED--semigroups makes easy to study the properties of
their $h$--gaps.
\begin{corollary}\label{cor5}{\rm $\;$}\\\vspace{-.3cm}
\begin{eqnarray}
&&\min\Delta_{{\cal H}}\left({\bf d}^m_{MED}\right)=d_m-d_{m-1}\;,\;\;\;\max
\Delta_{{\cal H}}\left({\bf d}^m_{MED}\right)=d_{m-1}-m\;,\label{ar25}\\
&&\#\Delta_{{\cal H}}\left({\bf d}^m_{MED}\right)=\frac{2}{m}\sum_{k=2}^md_k-
d_m\;,\;\;\;\#\Delta_{{\cal G}}\left({\bf d}^m_{MED}\right)=d_m-\frac{1}{m}
\sum_{k=2}^md_k-\frac{m-1}{2}\;.\;\;\;\;\label{ar25r}
\end{eqnarray}
\end{corollary}
{\sf Proof} $\;\;\;$Proof of (\ref{ar25}) follows from Lemma \ref{lem1} and
Theorem \ref{the7}. Proof of (\ref{ar25r}) follows from formulas (\ref{tr8}) 
and (\ref{tr9}) for cardinalities of $\Delta_{{\cal H}}\left({\bf d}^m\right)$
and $\Delta_{{\cal G}}\left({\bf d}^m\right)$ sets and formulas (\ref{ar24}) 
for the Frobenius number and the genus of the MED--semigroups.
$\;\;\;\;\;\;\Box$

\vspace{.1cm}
Formulas (\ref{ar25}) and (\ref{ar25r}) allow us to formulate two more
Corollaries.
\begin{corollary}\label{cor6}{\rm $\;$}\\
Let a numerical semigroup ${\sf S}\left({\bf d}^m_{MED}\right)$ be given. Then
\begin{eqnarray}
d_m\leq 2d_{m-1}-2m+3\;,\label{ar27}
\end{eqnarray}
and equality in (\ref{ar27}) holds iff $\#\Delta_{{\cal H}}\left({\bf d}^m_
{MED}\right)=1$.
\end{corollary}
{\sf Proof} $\;\;\;$By Theorem \ref{the8} the numerical semigroup ${\sf S}
\left({\bf d}^m_{MED}\right)$ is nonsymmetric, and therefore, by Corollary
\ref{cor2a}, $\max\Delta_{{\cal H}}\left({\bf d}^m_{MED}\right)\geq\min\Delta_
{{\cal H}}\left({\bf d}^m_{MED}\right)+t\left({\bf d}^m\right)-2$.
Substituting (\ref{ar24y}) and (\ref{ar25}) in the last inequality we get
(\ref{ar27}). Due to the equivalence of equalities
\begin{eqnarray}
\max\Delta_{{\cal H}}\left({\bf d}^m\right)=\min\Delta_{{\cal H}}\left({\bf 
d}^m\right)\;\;\;\Longleftrightarrow\;\;\;\#\Delta_{{\cal H}}\left({\bf 
d}^m\right)=1\;,\nonumber
\end{eqnarray}
we conclude that the equality in (\ref{ar27}) is attained iff $\#\Delta_{
{\cal H}}\left({\bf d}^m_{MED}\right)=1$.$\;\;\;\;\;\;\Box$
\begin{remark}
{\rm Note that the equality in (\ref{ar27}) can be fulfilled by
nonpseudo--symmetric MED--semigroup for $m\geq 4$ (see Remark \ref{re1} in
Section \ref{s1p5}). The inequality (\ref{ar27}) coincides with inequality
(\ref{ar24a}) in the case $m=3$.}
\end{remark}
The other Corollary deals with the cardinalities of sets $\Delta_{{\cal H}}
\left({\bf d}^m_{MED}\right)$ and $\Delta_{{\cal G}}\left({\bf d}^m_{MED}
\right)$.
\begin{corollary}\label{cor7}{\rm $\;$}\\
Let a numerical semigroup ${\sf S}\left({\bf d}^m_{MED}\right)$ be given. 
Then
\begin{eqnarray}
\#\Delta_{{\cal H}}\left({\bf d}^m_{MED}\right)=\#\Delta_{{\cal G}}\left({\bf 
d}^m_{MED}\right)\;,\;\;\;\;\mbox{iff}\label{ar26a}\\
6\sum_{k=2}^{m-1}d_k+m(m-1)=2(2m-3)d_m\;.\label{ar26}
\end{eqnarray}
\end{corollary}
{\sf Proof} $\;\;\;$By Theorem \ref{the8} the numerical semigroup ${\sf S}
\left({\bf d}^m_{MED}\right)$ is nonsymmetric, and therefore $\Delta_{{\cal 
H}}\left({\bf d}^m_{MED}\right)$ is nonempty. Equating both cardinalities 
$\#\Delta_{{\cal H}}\left({\bf d}^m_{MED}\right)$ and $\#\Delta_{{\cal G}}
\left({\bf d}^m_{MED}\right)$ given in (\ref{ar25r}) we get (\ref{ar26}). Vice 
versa, if (\ref{ar26}) does hold for the MED--semigroup then by comparison 
with (\ref{ar25r}) we get $\#\Delta_{{\cal H}}\left({\bf d}^m_{MED}\right)=
\#\Delta_{{\cal G}}\left({\bf d}^m_{MED}\right)$.$\;\;\;\;\;\;\Box$
\begin{remark}
{\rm By consequence of (\ref{ar26}), in the case $m=3$ equality (\ref{ar26a})
occurs iff
\begin{eqnarray}
{\bf d}^3_{MED}=(3,3k+1,3k+2)\;,\;\;\;k\geq 1\;.\label{ar28}
\end{eqnarray}}
\end{remark}
We finish this Section with one Example for the 5-dim numerical MED--semigroup.
The polynomial $Q\left({\bf d}^m;z\right)$ was obtained by means of 
diagrammatic calculation on the set $\Delta\left({\bf d}^m\right)$.
\begin{example}$\{d_1,d_2,d_3,d_4,d_5\}=\{5,7,9,11,13\}$\label{ex3}
{\footnotesize
\begin{eqnarray}
\left\{\begin{array}{rcr}
a_{11}d_1&=&20\\a_{22}d_2&=&14\\a_{33}d_3&=&18\\a_{44}d_4&=&22\\a_{55}d_5&=&
26\end{array}\right.,\;\;\;\;\;\;\left\{\begin{array}{lrl}F(5,7,9,11,13)&=&8\\
G(5,7,9,11,13)&=&6\\\Delta_{{\cal G}}(5,7,9,11,13)&=&\{1,3,8\}\\\Delta_{{\cal
H}}(5,7,9,11,13)
&=&\{2,4,6\}\\\#\left\{Q(5,7,9,11,13;z)\right\}&=&50\end{array}\right.,\;\;\;\;
\;\;\;\left.\begin{array}{lrl}\deg Q&=&53\\\max\Xi&=&51\\\min\Xi&=&14\\
\min\Delta_{{\cal H}}&=&53-51\end{array}\right.\;,\nonumber
\end{eqnarray}
\begin{eqnarray}
{\sf S}^{\prime}(5,7,9,11,13)&=&\{2,4,6,8\}\;,\;\;\;\;
\beta_1=10\;,\;\;\beta_2=20\;,\;\;\beta_3=15\;,\;\;\beta_4=4\;,\nonumber\\
H(5,7,9,11,13;z)&=&\frac{Q(5,7,9,11,13;z)}{(1-z^{5})
(1-z^{7})(1-z^{9})(1-z^{11})(1-z^{13})}\;,\;\;\;\;\Sigma_5=45\;,\nonumber\\
Q(5,7,9,11,13;z)&=&1-z^{14}-z^{16}-2z^{18}-2z^{20}-2z^{22}+z^{23}-
z^{24}+2z^{25}-z^{26}+3z^{27}+4z^{29}+4z^{31}+\nonumber\\
&&3z^{33}-z^{34}+2z^{35}-2z^{36}+z^{37}-3z^{38}-3z^{40}-3z^{42}-2z^{44}-
z^{46}+\nonumber\\
&&z^{47}+z^{49} + z^{51} + z^{53}\;.\nonumber
\end{eqnarray}}
\end{example}
\subsection{Semigroups of Maximal and Almost Maximal Length}\label{s3p5}
In this Section we apply our results of Section \ref{p4} to another special
kind of numerical semigroups, making use of important Theorem \cite{frob87}
on the structure of numerical semigroups. Recall Remark \ref{re0} on the
number of  nongaps of a semigroup.
\begin{theorem}{\rm (\cite{frob87}, Theorem 20)}\label{the9}\\
Let ${\sf S}\left({\bf d}^m\right)$ be a semigroup of type $t\left({\bf
d}^m\right)$ and let
$\#\Delta_{{\cal G}}\left({\bf d}^m\right)$ be the number of elements in
${\sf S}\left({\bf d}^m\right)$ which are less than $F\left({\bf d}^m\right)$.
Then we have
\begin{eqnarray}
G\left({\bf d}^m\right)\leq\#\Delta_{{\cal G}}\left({\bf d}^m\right)\cdot
t\left({\bf d}^m\right)\;.\label{bu1}
\end{eqnarray}
\end{theorem}
Notice that (\ref{bu1}) together with (\ref{tr8}) and (\ref{tr9}) is equivalent
to
\begin{eqnarray}
0\leq\#\Delta_{{\cal H}}\left({\bf d}^m\right)\leq\#\Delta_{{\cal G}}\left({\bf
d}^m\right)\cdot\left[t\left({\bf d}^m\right)-1\right]\;.\label{intro7}
\end{eqnarray}
Call the two extreme cases in (\ref{intro7}) {\em semigroups of maximal length}
(ML), when $\#\Delta_{{\cal H}}\left({\bf d}^m\right)$ attains its maximal 
value, and {\em semigroups of minimal length} ({\sl ml}), when $\#\Delta_{
{\cal H}}\left({\bf d}^m\right)$ vanishes. The latter case corresponds exactly 
to symmetric semigroups.
\begin{theorem}{\rm (\cite{brown91}, Theorem 1; \cite{brown92}, Corollary at
p. 339)}\label{the10}\\
Let ${\sf S}\left({\bf d}^m\right)$ be a numerical semigroup of type $t\left(
{\bf d}^m\right)$. Then ${\sf S}\left({\bf d}^m_{ML}\right)$ has maximal 
length iff
\begin{eqnarray}
{\bf d}^m_{ML}=\{m,mk+1,\ldots,mk+m-1\}\;,\;\;k\geq 1\;.\label{intro8a}
\end{eqnarray}
\end{theorem}
When $\#\Delta_{{\cal H}}\left({\bf d}^m\right)$ gets an intermediate values in 
the interval $\left(0,\;\#\Delta_{{\cal G}}\left({\bf d}^m\right)\cdot\left[t
\left({\bf d}^m\right)-1\right]\right)$, semigroup ${\sf S}\left({\bf d}^m
\right)$ can be generated by sophisticated series of generators, similar 
to (\ref{intro8a}), 
or by a significant number of sporadic tuples \cite{brown91}. In the meantime, 
the study of numerical semigroups with generic length are far from its 
completeness.

Call ${\sf S}\left({\bf d}^m\right)$ {\em semigroup of almost minimal length}
({\sl aml}) or {\em almost maximal length} (AML) if it has, respectively,
\begin{eqnarray}
\#\Delta_{{\cal H}}\left({\bf d}^m_{{\sl aml}}\right)=1\;,\;\;\;\;\mbox{or}\;
\;\;\;\#\Delta_{{\cal H}}\left({\bf d}^m_{AML}\right)=\#\Delta_{{\cal G}}
\left({\bf d}^m_{AML}\right)\cdot\left[t({\bf d}^m_{AML})-1\right]-1\;.
\label{intro8b}
\end{eqnarray}
The semigroups of almost minimal length were discussed in Section
\ref{s1p5}. As for the semigroups of almost maximal length,
several explicit results are known:
\begin{lemma}{\rm (\cite{brown92}, Proposition at p. 345)}\label{lem4}\\
Let ${\sf S}\left({\bf d}^m_{AML}\right)$ be a numerical semigroup of almost
maximal length and  $t\left({\bf d}^m_{AML}\right)$ its type. Then either
\begin{eqnarray}
i)\;\;t\left({\bf d}^m_{AML}\right)=2\;,\;\;\mbox{and}\;\;d_1=3\;\;\mbox{or}\;
\;4\;;\;\;\;\;\;\;ii)\;\;t\left({\bf d}^m_{AML}\right)\geq 3\;,\;\;\mbox{and}\;
\;d_1=t\left({\bf d}^m_{AML}\right)+1\;.\label{intro8c}
\end{eqnarray}
\end{lemma}
\begin{theorem}{\rm (\cite{brown91}, Theorem 3)}\label{the11}\\
Let ${\sf S}\left({\bf d}^m\right)$ be a numerical semigroup of type $t\left(
{\bf d}^m\right)=2$. Then ${\sf S}\left({\bf d}^m_{AML}\right)$ has almost
maximal length iff its generating tuple ${\bf d}^m_{AML}$ is one of the
following:
\begin{eqnarray}
i)\;{\bf d}^m_{AML}=\{3,3k+2,3k+4\},\;k\geq 1\;,\;\;ii)\;{\bf d}^m_{AML}=
\{4,5,11\},\;\;\;iii)\;{\bf d}^m_{AML}=\{4,7,13\}.\label{intro9a}
\end{eqnarray}
\end{theorem}
\begin{theorem}{\rm (\cite{brown91}, Theorem 5; \cite{brown92}, Corollary 2)}
\label{the12}\\
Let ${\sf S}\left({\bf d}^m\right)$ be a numerical semigroup of type $t\left(
{\bf d}^m\right)=t_m\geq 3$. Then ${\sf S}\left({\bf d}^m_{AML}\right)$
has almost maximal length iff
\begin{eqnarray}
{\bf d}^m_{AML}=\left\{t_m+1,k(t_m+1)+t_m,k(t_m+1)+(t_m+2),k(t_m+1)+(t_m+3),
\ldots,k(t_m+1)+2t_m\right\}.\nonumber
\end{eqnarray}
\end{theorem}
We obtain two Corollaries on the semigroups of maximal and almost maximal
length respectively.
\begin{corollary}\label{cor8}{\rm $\;$}\\
The numerical semigroups ${\sf S}\left({\bf d}^m_{ML}\right)$ satisfy
\begin{eqnarray}
\left\{\begin{array}{l}F\left({\bf d}^m_{ML}\right)=mk-1\\\#\Delta_{{\cal H}}
\left({\bf d}^m_{ML}\right)=mk-2k\\\#\Delta_{{\cal G}}\left({\bf d}^m_{ML}
\right)=k\end{array}\right.,\;\;\;\;\;\left\{\begin{array}{l}
\min\Delta_{{\cal H}}\left({\bf d}^m_{ML}\right)=1\\\max\Delta_{{\cal H}}
\left({\bf d}^m_{ML}\right)=mk-2\end{array}\right.\;.\label{intro10}
\end{eqnarray}
\end{corollary}
{\sf Proof} $\;\;\;$By Theorem \ref{the10} any ML--semigroup is a
MED--semigroup; therefore, according to (\ref{ar24}) and (\ref{ar25}), we
arrive at (\ref{intro10}).$\;\;\;\;\;\;\Box$
\begin{remark}
{\rm By consequence of (\ref{intro10}), there exists only one
ML--semigroup with $\#\Delta_{{\cal H}}\left({\bf d}^m_{ML}\right)=1$
which is pseudo--symmetric, ${\sf S}(3,4,5 )$.}
\end{remark}

Note that the AML--semigroup ${\sf S}(3,3k+2,3k+4)$ in Theorem \ref{the11}
is also a MED--semigroup,
\begin{eqnarray}
\left\{\begin{array}{l}F(3,3k+2,3k+4)=3k+1\\
\#\Delta_{{\cal H}}(3,3k+2,3k+4)=k\\\#\Delta_{{\cal G}}
(3,3k+2,3k+4)=k+1\end{array}\right.,\;\;\;\;
\left\{\begin{array}{l}\min\Delta_{{\cal H}}(3,3k+2,3k+4)=2\\
\max\Delta_{{\cal H}}(3,3k+2,3k+4)=3k-1\end{array}\right.\;.
\label{intro11}
\end{eqnarray}
while the other sporadic AML--semigroups (\ref{intro9a}) with $t\left({\bf
d}^m_{AML}\right)=2$ are not MED--semigroups,
\begin{eqnarray}
\left\{\begin{array}{l}F(4,5,11)=7\\
\#\Delta_{{\cal H}}(4,5,11)=2\\\#\Delta_{{\cal G}}(4,5,11
)=3\end{array}\right.\left.\begin{array}{l}\min\Delta_{{\cal H}}
(4,5,11)=1\\\max\Delta_{{\cal H}}
(4,5,11)=6\end{array}\right.\;;\;\left\{\begin{array}{l}
F(4,7,13)=10\\
\#\Delta_{{\cal H}}(4,7,13)=3\\\#\Delta_{{\cal G}}(4,7,13
)=4\end{array}\right.\left.
\begin{array}{l}\min\Delta_{{\cal H}}(4,5,11)=1\\\max\Delta_
{{\cal H}}(4,7,13)=9\end{array}\right.\nonumber
\end{eqnarray}
\begin{corollary}\label{cor9}{\rm $\;$}\\
Any semigroup of almost maximal length ${\sf S}\left({\bf d}^m_{AML}\right)$
with $t_m>2$ satisfies
\begin{eqnarray}
\left\{\begin{array}{l}F\left({\bf d}^m_{AML}\right)=k(t_m+1)+t_m-1\\
\#\Delta_{{\cal H}}\left({\bf d}^m_{AML}\right)=(k+1)t_m-k-2\\
\#\Delta_{{\cal G}}\left({\bf d}^m_{AML}\right)=k+1\end{array}\right.,\;\;
\left\{\begin{array}{l}\min\Delta_{{\cal H}}\left({\bf d}^m_{AML}\right)=
1\\\max\Delta_{{\cal H}}\left({\bf d}^m_{AML}\right)=k(t_m+1)+t_m-2\end{array}
\right.\label{intro12}
\end{eqnarray}
\end{corollary}
{\sf Proof} $\;\;\;$
By Theorem \ref{the12} the ${\sf S}\left({\bf d}^m_{AML}\right)$--semigroup
with $t_m>2$ is the MED--semigroup, and therefore, according to (\ref{ar24})
and (\ref{ar25}), we arrive at (\ref{intro12}).$\;\;\;\;\;\;\Box$

\vspace{.1cm}
\begin{remark}
{\rm By consequence of (\ref{intro11}) and  (\ref{intro12}),
we get the known result \cite{baru97}: there exists only one AML--
pseudo--symmetric semigroup, ${\sf S}(3,5,7)$.}
\end{remark}
\section*{Acknowledgement}
The paper has been written during the stay of one of the authors (LGF) at the
Mathematics Section of the Abdus Salam International Centre for Theoretical
Physics (ICTP) and its hospitality is highly appreciated. The other author (FA)
is grateful to V. I. Arnold for having introduced her to the Frobenius problem
during his stay at ICTP.


\vspace{.3cm}
\begin{thebibliography}{99}
\bibitem{baru97} V. Barucci, D. E. Dobbs and M. Fontana, {\it Maximality
Properties in Numerical Semigroups and Applications to
One--Dimensional Analytically Irreducible Local Domains}, \\
Memoirs Amer. Math. Soc., {\bf 125, (598)}, 1997
\bibitem{abhy67}  S. Abhyankar, {\it Local Rings of High Embedding Dimension},
                   \\ Amer. J. Math, {\bf 89}, 1073 (1967)
\bibitem{stur05}  E. Miller and B. Sturmfels, {\it Combinatorial Commutative
                  Algebra}, \\ New York : Springer, 2005, p. 45
\bibitem{frob87}  R. Fr\"oberg, C. Gottlieb and R. H\"aggkvist, {\it On
                  Numerical Semigroups}, \\Semigroup Forum, {\bf 35}, 63 (1987)
\bibitem{brown91} W. C. Brown and F. Curtis, {\it Numerical Semigroups of
                  Maximal and Almost Maximal Length}, \\
                  Semigroup Forum, {\bf 42}, 219 (1991)
\bibitem{brown92} W. C. Brown and J. Herzog, {\it 1-dim Local Rings of
                  Maximal and Almost Maximal Length}, \\
                  J. Algebra, {\bf 151}, 332 (1992)
\bibitem{stan96}  R. P. Stanley, {\it Combinatorics and Commutative Algebra},\\
                  Birkh\"auser Boston, 1996, 2nd ed.
\bibitem{heku71}  J. Herzog and E. Kunz, {\it Die Werthalbgruppe eines lokalen
                  Rings der Dimension 1},\\ Sitzungsberichte der Heidelberger
                  Akademie der Wissenschaften, Springer, Berlin (1971)
\bibitem{aper46}  R. Ap\'ery, {\em Sur les Branches Superlin\'eaires des
                  Courbes Alg\'ebriques},\\ C. R. Acad. Sci. Paris,
                  {\bf 222}, 1198 (1946). MR {\bf 8}, 221
\bibitem{herz70}  J. Herzog, {\it Generators and Relations of Abelian
                  Semigroups and Semigroup Rings},\\
                  Manuscripta Math., {\bf 3}, 175 (1970)
\bibitem{frob94}  R. Fr\"oberg, {\it The Frobenius Number of Some Semigroups},
                  \\ Comm. Algebra, {\bf 22} (14), 6021 (1994)
\bibitem{denh03}  G. Denham, {\it Short Generating Functions for Some
                  Semigroup Algebras}, \\Electronic J. Combinatorics,
                  {\bf 10} (2003), \# R36, 7 pages.
\bibitem{fel04}   L. G. Fel, {\it Frobenius Problem for Semigroups ${\rm S}
                  (d_1,d_2,d_3)$},\\ Funct. Analysis and Other Math.,
                  {\bf 1}, \# 2, 135 (2006), \\
                  preprint, http://arxiv.org/abs/math.NT/0409331
\bibitem{john60}  S. M. Johnson, {\it A Linear Diophantine Problem}, \\
                  Canad. J. Math., {\bf 12}, 390 (1960)
\bibitem{worm86}  L. A. Sz\'ekely and N. C. Wormald, {\it Generating Functions
                  for the Frobenius Problem \\with 2 and 3 Generators}, \\
                  Math. Chronicle, {\bf 15}, 49 (1986)
\bibitem{wag05}   D. Beihoffer, J. Hendry, A. Nijenhuis and S. Wagon, {\em
                  Faster Algorithms for Frobenius \\Numbers}, \\ Electronic
                  J. of Combinatorics, {\bf 12}, R27 (2005)
\bibitem{ber79}   H. Bresinsky, {\it Monomial Gorenstein Ideals}, \\
                  Manuscripta Math., {\bf 29}, 159 (1979)
\bibitem{ros00}   J. C. Rosales, {\it Numerical Semigroups with Ap\'ery Sets
                  of Unique Expression}, \\J. Algebra, {\bf 226}, 479 (2000)
\bibitem{kill00}  H. G. Killingbergtr\O, {\em Betjening av figur i Frobenius'
                  problem}, (In Norwegian), \\ Normat (Nordisk Matematisk
                  Tidskrift), {\bf 2}, 75 (2000)
\bibitem{sal79}   J. D. Sally, {\it Cohen--Macaulay Local Rings of Maximal
                  Embedding Dimensions}, \\ J. Algebra, {\bf 56}, 168 (1979)
\end{thebibliography}
\end{document}